\documentclass[11pt,letterpaper]{amsart}

\usepackage{mathrsfs,amssymb,esint,verbatim,microtype}
\numberwithin{equation}{section}

\newtheorem{thm}{Theorem}[section]
\newtheorem{cor}[thm]{Corollary}
\newtheorem{lem}[thm]{Lemma}
\newtheorem{prop}[thm]{Proposition}

\usepackage[pagebackref=true,colorlinks=true,bookmarksopen,bookmarksnumbered,
citecolor=red,linkcolor=blue,urlcolor=cyan]{hyperref}
\hypersetup{
  pdftitle={Dirichlet Green kernel estimates and Sobolev-type inequalities for twisted differential forms},
  pdfauthor={Fusheng Deng, Gang Huang, and Xiangsen Qin},
  pdfsubject={Dirichlet Green kernels for twisted Dolbeault and de Rham complexes},
  pdfkeywords={Dirichlet Green kernel, Dolbeault Laplacian, Sobolev-type inequalities}
}
\allowdisplaybreaks[4]
\begin{document}

\title[Dirichlet Green kernel estimates]
{Dirichlet Green kernel estimates and Sobolev-type inequalities for twisted differential forms}

\author[F. Deng]{Fusheng Deng}
\address{School of Mathematical Sciences, University of Chinese Academy of Sciences \\ Beijing 100049, P. R. China}
\email{fshdeng@ucas.ac.cn}

\author[G. Huang]{Gang Huang}
\address{School of Mathematical Sciences, University of Chinese Academy of Sciences \\ Beijing 100049, P. R. China}
\email{huanggang21@mails.ucas.ac.cn}

\author[X. Qin]{Xiangsen Qin}
\address{Chern Institute of Mathematics and LPMC, Nankai University \\ Tianjin 300071, China}
\email{qinxiangsen@nankai.edu.cn}

\begin{abstract}
We study the full-trace Dirichlet realization of the Dolbeault Laplacian on differential forms with values in a
Hermitian holomorphic vector bundle   over a relatively compact smooth domain in
a K\"ahler manifold. We prove global Green kernel estimates that are uniform up to the boundary, including one- and two-boundary-factor bounds and estimates for the \(\bar\partial\)- and \(\bar\partial^*\)-derivatives. These estimates yield Sobolev-type inequalities with boundary terms. In real dimension two, the first-order bound has a logarithmic loss. In top antiholomorphic degree, we obtain quantitative \(L^r\)-to-\(L^k\) solvability for \(\bar\partial\) without pseudoconvexity. Analogous  results hold for
the twisted de Rham complex of a flat metric connection on a compact
Riemannian manifold with smooth boundary.
\end{abstract}

\subjclass[2020]{Primary 35J08, 32W05; Secondary 58J32, 58J35, 46E35}
\keywords{Dirichlet Green kernel, Dolbeault Laplacian, Sobolev-type inequalities}

\maketitle
\tableofcontents
\section{Introduction}

Dirichlet Green kernels reflect both the singular behavior of an elliptic inverse near the diagonal and the decay imposed by the Dirichlet boundary condition. For scalar uniformly elliptic equations, Gr\"uter--Widman established
Green-function estimates, Davies related Green-function and heat-kernel
bounds, and Zhang obtained sharp short-time Dirichlet heat-kernel estimates
with boundary-distance factors on bounded \(C^{1,1}\) domains
\cite{GW82,D87,D89,Z02}. For elliptic systems, classical Green-matrix
results include \cite{S70,F86,DM95}; subsequent work treated rougher
coefficients and domains under appropriate regularity hypotheses
\cite{HK07,DK09,KK10}, while Conlon--Giunti--Otto proved almost-everywhere
existence for measurable systems and averaged pointwise estimates in the
stationary setting \cite{CGO17}. Li--Strohmaier obtained interior
“not-feeling-the-boundary” estimates for arbitrary nonnegative self-adjoint
extensions of the componentwise Euclidean Laplacian \cite{LS16}, but these
do not yield Dirichlet boundary factors. A recent preprint of
Hsiao--Marinescu--Zhu studies semiclassical heat-kernel asymptotics near the
boundary for the \(\bar\partial\)-Neumann Laplacian in high tensor powers
under condition \(Z(q)\) \cite{HMZ26}; its realization and asymptotic regime differ
from the fixed-bundle full-trace Dirichlet problem considered here. General
smooth boundary realizations and Hodge boundary problems are developed in
\cite{G96,S95,MMT01,MMMT16}.

For differential forms, the word ``Dirichlet'' is not unambiguous. In this
paper it always refers to the full-trace realization
\[
\operatorname{Dom}(L)=H^2\cap H_0^1,
\]
so every coefficient of the form has zero boundary trace. This is a strongly
elliptic realization, but it is neither the \(\bar\partial\)-Neumann
realization nor the absolute or relative realization used in Hodge theory.
The distinction matters. The latter boundary conditions are adapted to the
underlying differential complex; the full-trace condition is not. In return,
it has the scalar Dirichlet boundary behavior that is central to the estimates
proved here.

Our first result gives global Green kernel bounds, uniform up to the boundary,
for the full-trace Dirichlet realization of the Dolbeault Laplacian.
Derivatives in the second variable act on the dual bundle factor; the precise
convention is given before Theorem~\ref{thm:integral representations}.

\begin{thm}\label{thm: green form complex boundary}
Let $E$ be a Hermitian holomorphic vector bundle over a K\"ahler manifold
$(M,g)$ of complex dimension $n$, let $\Omega\Subset M$ be a domain with
nonempty smooth boundary, and fix $p,q\in\{0,\ldots,n\}$. Let
$\square_{p,q}$ be the $\bar\partial$-Laplacian on $E$-valued $(p,q)$-forms,
and denote its Dirichlet realization on
$L^2(\Omega,\Lambda^{p,q}T^*M\otimes E)$ by the same symbol. Thus
\[
\operatorname{Dom}(\square_{p,q})
=H^2(\Omega,\Lambda^{p,q}T^*M\otimes E)
\cap H_0^1(\Omega,\Lambda^{p,q}T^*M\otimes E).
\]
Put $\delta(z):=d(z,\partial\Omega)$, set
$\mathcal G_{p,q}:=\square_{p,q}^{-1}$, and denote the Schwartz kernel of
$\mathcal G_{p,q}$ by $G_{p,q}(x,y)$. Then there is a constant
\(C=C(\overline\Omega,g,E)>0\) such that, for all distinct
\(x,y\in\overline\Omega\), the following estimates hold:
\begin{itemize}
\item[(i)]
\[
|G_{p,q}(x,y)|\leq C\begin{cases}
d(x,y)^{2-2n}, & n\geq2,\\
1+|\ln d(x,y)|, & n=1;
\end{cases}
\]
\item[(ii)]
\begin{align*}
|G_{p,q}(x,y)|
&\leq C\min\{\delta(x),\delta(y)\}d(x,y)^{1-2n},\\
|G_{p,q}(x,y)|
&\leq C\delta(x)\delta(y)d(x,y)^{-2n};
\end{align*}
\item[(iii)]
{\small
\[
\max\!\left\{|\bar\partial_yG_{p,q}(x,y)|,
|\bar\partial_y^{*}G_{p,q}(x,y)|\right\}
\leq C\begin{cases}
d(x,y)^{1-2n},&n\geq2,\\
\bigl(1+|\ln d(x,y)|\bigr)d(x,y)^{-1},&n=1.
\end{cases}
\]
}
and
\[
\max\!\left\{|\bar\partial_yG_{p,q}(x,y)|,
|\bar\partial_y^{*}G_{p,q}(x,y)|\right\}
\leq C\delta(x)d(x,y)^{-2n}.
\]
\end{itemize}
\end{thm}

The zero-order singularities in Theorem~\ref{thm: green form complex boundary}
have the order of the Euclidean fundamental solution in real dimension
\(2n\). Part~(ii) records both the one-boundary-factor and the
two-boundary-factor estimates. The weighted estimate in part~(iii) has Poisson
scale when the differentiated variable approaches the boundary. In complex
dimension one, our argument introduces a logarithmic loss in the unweighted
first-order estimate. We do not know whether that loss is necessary.

The scalar version of these bounds was proved in \cite{DHQ24}, while
\cite{DHQ25} treats twisted differential forms on compact K\"ahler manifolds
without boundary.

The proof separates short and long times in the heat representation of the
Green operator. For short times, the scalar Dirichlet heat kernel supplies the
boundary factors, and semigroup domination transfers the estimate to the
bundle-valued heat kernel \cite{HSU77}. For long times, we use the Dirichlet
eigenform expansion. The heat-trace argument of Donnelly--Li gives the required
eigenvalue growth \cite{DL82}; the boundary-weighted eigenform bounds are then
obtained by combining semigroup domination with the scalar boundary estimate.
The full-trace kernel is trivial, by unique continuation for a Dirac-type
operator, so the first Dirichlet eigenvalue is positive and the long-time part
decays exponentially. This short-time/long-time separation is why an arbitrary
fixed lower bound for the Weitzenb\"ock curvature is sufficient. No
nonnegativity assumption is needed.

The differentiated Green kernel is obtained from a uniform local estimate. If
\(\square_{p,q}f=0\) on a geodesic ball \(B(x_0,r)\Subset\Omega\), then
\[
\max\left\{
\sup_{B(x_0,r/2)}|\bar\partial f|,
\sup_{B(x_0,r/2)}|\bar\partial^*f|
\right\}
\leq \frac{C}{r}\sup_{B(x_0,r)}|f|.
\]
The proof adapts the Moser iteration of Lu, Qi S. Zhang, and Zhu
\cite[Proposition~3.3]{LZZ21} to bundle-valued forms on uniformly controlled
geodesic balls. A Euclidean-domain version, with fully tracked constants,
appears in \cite[Lemma~5.5]{Q25}. Section~\ref{sec:green form complex estimate}
keeps the constants explicit enough to display their dependence on the 
manifold geometry, a common lower bound of the Weitzenb\"ock curvature, and the first Dirichlet eigenvalue of $\square_{p,q}$.

Green representation turns the kernel estimates into Sobolev inequalities
with boundary traces. The volume kernels have Riesz orders two and one in real
dimension \(2n\), while the boundary term is controlled by the Poisson-scale
kernel
\[
(x,y)\longmapsto \frac{\delta(x)}{d(x,y)^{2n}},
\qquad x\in\Omega,\quad y\in\partial\Omega.
\]
The exact bundle-valued representation formulas, including the dual-variable
conventions and the boundary operator, are given in
Section~\ref{section:real inequality}.

To state the resulting ranges, let \(m\geq2\), \(0<\alpha\leq m\), and
\(1\leq a,b\leq\infty\). We call \((a,b)\) \((m,\alpha)\)-admissible if
\begin{equation}\label{equ:admissible-volume}
\begin{cases}
1\leq b<\dfrac{m}{m-\alpha},&a=1,\ \alpha<m,\\[4pt]
1\leq b<\dfrac{2a}{2-a},
&(m,\alpha)=(2,1),\ 1<a<2,\\[4pt]
1\leq b\leq\dfrac{ma}{m-\alpha a},
&(m,\alpha)\neq(2,1),\ 1<a<\dfrac m\alpha,\\[4pt]
1\leq b<\infty,&a=\dfrac m\alpha,\\[4pt]
1\leq b\leq\infty,&a>\dfrac m\alpha.
\end{cases}
\end{equation}
We call \((a,b)\) \(m\)-boundary-admissible if
\begin{equation}\label{equ:admissible-boundary}
\begin{cases}
1\leq b<\dfrac{m}{m-1},&a=1,\\[4pt]
1\leq b\leq\dfrac{ma}{m-1},&1<a\leq\infty.
\end{cases}
\end{equation}
For the two-dimensional first-order volume kernel, the logarithmic loss in
Theorem~\ref{thm: green form complex boundary}(iii) makes the critical
strong-type range strict. It does not change the boundary range, because the
boundary kernel has the Poisson-scale bound without a logarithmic factor.

\begin{cor}\label{cor: complex laplace}
Under the notation of Theorem~\ref{thm: green form complex boundary}, let
$1\leq k,r,s\leq\infty$. Assume that $(r,k)$ is $(2n,2)$-admissible and
$(s,k)$ is $2n$-boundary-admissible. Then there is a constant
$C=C(\overline\Omega,g,E,k,r,s)>0$ such that
\[
\|f\|_{L^k(\Omega)}
\leq C\left(
\|\square_{p,q}f\|_{L^r(\Omega)}
+\|f\|_{L^s(\partial\Omega)}\right)
\]
for every
$f\in C^2(\overline\Omega,\Lambda^{p,q}T^*M\otimes E)$.
\end{cor}

\begin{cor}\label{cor: complex dbar}
Under the notation of Theorem~\ref{thm: green form complex boundary}, let
$1\leq k,r,s,t\leq\infty$. Assume that both $(r,k)$ and $(s,k)$ are
$(2n,1)$-admissible and that $(t,k)$ is $2n$-boundary-admissible. Then there is
a constant $C=C(\overline\Omega,g,E,k,r,s,t)>0$ such that
\[
\|f\|_{L^k(\Omega)}
\leq C\left(
\|\bar\partial f\|_{L^r(\Omega)}
+\|\bar\partial^*f\|_{L^s(\Omega)}
+\|f\|_{L^t(\partial\Omega)}\right)
\]
for every
$f\in C^1(\overline\Omega,\Lambda^{p,q}T^*M\otimes E)$.
\end{cor}

The passage from a Green representation to Sobolev inequalities with boundary
traces is classical. Cianchi and Maz'ya established a substantially more
general Euclidean theory on arbitrary domains, with geometry-independent
constants and sharp target spaces \cite{CM16}. Shi and Yao recently constructed
a \(\bar\partial\)-solution operator that gains one half derivative in
fractional Sobolev spaces on bounded strictly pseudoconvex domains with
\(C^2\)-boundary \cite{SY25}. Duse and Ros\'en obtained coercive and Morrey
estimates for broad classes of first-order systems under natural half-boundary
conditions \cite{DR26}. These results concern different boundary realizations
or different analytic regimes. The specific input in Corollaries
\ref{cor: complex laplace} and \ref{cor: complex dbar} is the full-trace,
bundle-valued Green kernel of Theorem~\ref{thm: green form complex boundary},
including its two boundary factors and its differentiated Poisson-scale bound.

Top antiholomorphic degree gives a further consequence. Every \((p,n)\)-form
is automatically \(\bar\partial\)-closed, and the full-trace Green operator
therefore produces a right inverse for \(\bar\partial\) without a
pseudoconvexity assumption.

\begin{thm}\label{thm:complex L^p estimate}
Under the notation of Theorem~\ref{thm: green form complex boundary}, let
\(1\leq r,k\leq\infty\), and assume that \((r,k)\) is
\((2n,1)\)-admissible. Set $\mathcal S_{p,n}:=\bar\partial^*\mathcal G_{p,n}$. This operator has
a canonical bounded extension
\[
\mathcal S_{p,n}:L^r(\Omega,\Lambda^{p,n}T^*M\otimes E)
\longrightarrow
L^k(\Omega,\Lambda^{p,n-1}T^*M\otimes E),
\]
and
\begin{equation*}
\bar\partial(\mathcal S_{p,n}f)=f
\quad\text{in }\mathcal D'(\Omega,\Lambda^{p,n}T^*M\otimes E),
\qquad
\|\mathcal S_{p,n}f\|_{L^k(\Omega)}
\leq C_{r,k}\|f\|_{L^r(\Omega)},
\end{equation*}
where $C_{r,k}$ depends only on
$(\overline\Omega,g,E,r,k)$.
\end{thm}

For \(n\geq2\), the critical strong exponent
\(k=2nr/(2n-r)\) is attained when \(1<r<2n\). If \(r=2n\), every
finite \(k\) is allowed, and \(r>2n\) permits \(k=\infty\). In complex
dimension one, the logarithmic loss gives the strict ranges
\[
k<2\quad(r=1),
\qquad
k<\frac{2r}{2-r}\quad(1<r<2).
\]
For \(r=2\), every finite \(k\) is allowed, while \(r>2\) permits
\(k=\infty\).

The classical \(L^2\) theory on pseudoconvex domains was established in
foundational work of Kohn, H\"ormander, and Folland--Kohn
\cite{K63,H65,FK72}. On strongly pseudoconvex
domains, Kerzman, {\O}vrelid, Krantz, and Beals--Greiner--Stanton obtained
\(L^p\), H\"older, and Lipschitz estimates through integral formulas and the
\(\bar\partial\)-Neumann framework \cite{K71,O71,K76,BGS87}. Li treats
complete K\"ahler manifolds without boundary under assumptions on Weitzenb\"ock curvature \cite{LX10}. Theorem~\ref{thm:complex L^p estimate} is not
a substitute for these complex-adapted theories. It applies on every smooth
relatively compact K\"ahler domain, but only in top antiholomorphic degree. If
\(v=\mathcal G_{p,n}f\), then \(\bar\partial v=0\) for degree reasons. For
\(q<n\), the full-trace Green operator need not preserve
\(\ker\bar\partial\), so no analogous lower-degree solvability statement is
made here.

The same method applies to the twisted de Rham complex for a flat metric
connection. The real counterpart of
Theorem~\ref{thm: green form complex boundary} is as follows.

\begin{thm}\label{thm: green form boundary}
Let $(M,g)$ be a compact Riemannian manifold of dimension $n\geq2$ with
nonempty smooth boundary, let $E$ be a Hermitian vector bundle over $M$ with a
flat metric connection $D$, and fix $p\in\{0,\ldots,n\}$. Let $\Delta_p$ be
the Hodge Laplacian on $E$-valued $p$-forms, with full-trace Dirichlet domain
\[
\operatorname{Dom}(\Delta_p)
=H^2(M,\Lambda^pT^*M\otimes E)
\cap H_0^1(M,\Lambda^pT^*M\otimes E).
\]
Put $\delta(z):=d(z,\partial M)$, let
$\mathcal G_p:=\Delta_p^{-1}$, and denote its Schwartz kernel by $G_p(x,y)$.
Then there is a constant $C=C(M,g,E,D)>0$ such that, for all distinct
$x,y\in M$, the following estimates hold:
\begin{itemize}
\item[(i)]
\[
|G_p(x,y)|\leq C\begin{cases}
d(x,y)^{2-n}, & n\geq3,\\
1+|\ln d(x,y)|, & n=2;
\end{cases}
\]
\item[(ii)]
\begin{align*}
|G_p(x,y)|
&\leq C\min\{\delta(x),\delta(y)\}d(x,y)^{1-n},\\
|G_p(x,y)|
&\leq C\delta(x)\delta(y)d(x,y)^{-n};
\end{align*}
\item[(iii)]
\[
\max\{|d_yG_p(x,y)|,|d_y^*G_p(x,y)|\}
\leq C\begin{cases}
d(x,y)^{1-n}, & n\geq3,\\
\bigl(1+|\ln d(x,y)|\bigr)d(x,y)^{-1}, & n=2,
\end{cases}
\]
and
\[
\max\{|d_yG_p(x,y)|,|d_y^*G_p(x,y)|\}
\leq C\delta(x)d(x,y)^{-n}.
\]
\end{itemize}
\end{thm}

Scott developed an \(L^p\) theory for differential forms \cite{SC95}, and
Gol'dshtein--Troyanov related Sobolev inequalities to
\(L_{q,p}\)-cohomology \cite{GT06}. Wang proved Calder\'on--Zygmund estimates
for forms on closed manifolds \cite{W04}, while Li developed \(L^p\) Hodge
theory on complete manifolds under assumptions on Riesz transforms,
potentials, and Weitzenb\"ock curvature \cite{L09,L10}. More recently,
Gol'dshtein--Kopylov--Panenko estimated differential-form Sobolev embeddings on
Euclidean balls and their bi-Lipschitz images \cite{GKP26}. None of these
results gives the full-trace one- and two-boundary-factor estimates or the
Poisson-scale differentiated bound in
Theorem~\ref{thm: green form boundary}.

The proof is the same short-time/long-time argument as in the complex case,
with real dimension \(n\) in place of \(2n\). Flatness of $D$ allows $d$ and
$d^*$ to commute with $\Delta_p$, and the local Moser estimate gives
\[
\max\left\{
\sup_{B(x_0,r/2)}|df|,
\sup_{B(x_0,r/2)}|d^*f|
\right\}
\leq \frac{C}{r}\sup_{B(x_0,r)}|f|
\]
for \(\Delta_p\)-harmonic forms on interior balls. The corresponding Green
representation formulas and the same potential estimates yield the following
corollaries.

\begin{cor}\label{cor: real laplace}
Under the notation of Theorem~\ref{thm: green form boundary}, let
$1\leq q,k,r\leq\infty$, assume that $(k,q)$ is $(n,2)$-admissible, and assume
that $(r,q)$ is $n$-boundary-admissible. Then there is a constant
$C=C(M,g,E,D,q,k,r)>0$ such that
\[
\|f\|_{L^q(M)}
\leq C\left(
\|\Delta_p f\|_{L^k(M)}
+\|f\|_{L^r(\partial M)}\right)
\]
for every $f\in C^2(M,\Lambda^pT^*M\otimes E)$.
\end{cor}

\begin{cor}\label{cor: real gradient}
Under the notation of Theorem~\ref{thm: green form boundary}, let
$1\leq q,k,r,s\leq\infty$, assume that both $(k,q)$ and $(r,q)$ are
$(n,1)$-admissible, and assume that $(s,q)$ is
$n$-boundary-admissible. Then there is a constant $C>0$, depending only on
$(M,g,E,D,q,k,r,s)$, such that
\[
\|f\|_{L^q(M)}
\leq C\left(
\|df\|_{L^k(M)}
+\|d^*f\|_{L^r(M)}
+\|f\|_{L^s(\partial M)}\right)
\]
for every $f\in C^1(M,\Lambda^pT^*M\otimes E)$.
\end{cor}

Taking $q=r=\infty$ and $k>n/2$ in
Corollary~\ref{cor: real laplace} gives a boundary $L^\infty$ estimate.
Corollary~\ref{cor: real gradient} is used in \cite{HJQ25} to obtain improved
$L^2$ estimates for $d$. Related inequalities on complete manifolds without
boundary appear in \cite[Theorem~1.2]{WZ12} and
\cite[Theorem~2.1]{L09}.

The paper is organized as follows. Section~\ref{notation} fixes the geometric,
operator, and kernel conventions. Section~\ref{sec:green form complex estimate}
proves the local first-order estimate and the Dirichlet Green kernel bounds.
Section~\ref{section:real inequality} records the Green representation formulas,
develops the volume and boundary potential estimates, and proves
Corollaries~\ref{cor: complex laplace} and \ref{cor: complex dbar}.
Section~\ref{sec:L^p estimate} applies the differentiated kernel bound to the
top-degree \(\bar\partial\)-equation.

\subsection*{Acknowledgments}
This work was supported by the National Key R\&D Program of China
(No.~2021YFA1003100), the NSFC (No.~12471079), and the Fundamental Research
Funds for the Central Universities.

\section{Notation and Conventions}\label{notation}

We use \(\mathbb N=\{0,1,2,\ldots\}\), and all Hermitian inner products are
complex linear in the first argument. All manifolds are smooth and oriented;
their boundaries are smooth and carry the outward-normal-first orientation.
We denote the interior of \(M\) by \(M^\circ\), write \(A\Subset B\) if
\(\overline A\) is a compact subset of \(B\), and use \(dV\), \(dS\), and
\(d(x,y)\) for Riemannian volume, hypersurface measure, and distance. We set
\[
B(x,r):=\{y\in M:d(x,y)<r\}.
\]

For a metric vector bundle \(F\to M\), the notation
\[
C^k(U,F),\quad C_c^k(U,F),\quad L^p(U,F),\quad H^k(U,F)
\]
is standard; compact support is understood to lie in \(U^\circ\), and
\(C^k(\overline U,F)\) denotes sections extending to a neighborhood of
\(\overline U\). We suppress \(F\) from \(L^p\)-norms when it is clear and
write
\[
H_0^1(U,F):=
\overline{C_c^\infty(U^\circ,F)}^{\,H^1(U,F)}.
\]
For smooth boundary this is the subspace with zero trace. All formal adjoints
are taken with respect to
\[
(u,v)_{L^2(U)}:=\int_U\langle u,v\rangle\,dV.
\]
If \(T\) is self-adjoint on a finite-dimensional inner-product space, then
\(T\geq b\) means
\[
\langle Tu,u\rangle\geq b\langle u,u\rangle
\qquad\text{for all }u.
\]

Let \((M,g)\) have real dimension \(m\), and let \(E\to M\) be a metric
vector bundle with metric connection \(D\), also denoting the induced
connection on \(\Lambda^pT^*M\otimes E\). Let \(d_D\) be the exterior
covariant derivative and \(d_D^*\) its formal adjoint. Once \(D\) is fixed,
we write
\[
d:=d_D,\qquad d^*:=d_D^*,\qquad
\Delta_p:=dd^*+d^*d.
\]
The untwisted Hodge Laplacian is recovered by taking \(E\) to be the trivial
line bundle. If \(R^D\) is the curvature of the induced connection and
\(X_1,\ldots,X_m\) is a local orthonormal frame, define
\begin{equation}\label{equ:operator}
\begin{aligned}
(\mathfrak{Ric}_p^Ef)(Y_1,\ldots,Y_p)
:={}&
\sum_{i=1}^p\sum_{j=1}^m
(R^D(X_j,Y_i)f)\\
&\qquad
(Y_1,\ldots,Y_{i-1},X_j,Y_{i+1},\ldots,Y_p),
\end{aligned}
\end{equation}
with \(\mathfrak{Ric}_0^E=0\). Then
\[
\Delta_p=D^*D+\mathfrak{Ric}_p^E
\]
and
\[
\Delta|f|^2
=
-2\operatorname{Re}\langle\Delta_pf,f\rangle
+2|Df|^2
+2\operatorname{Re}\langle\mathfrak{Ric}_p^Ef,f\rangle,
\]
where \(D^*D\) is the nonnegative connection Laplacian and
\[
\Delta:=\operatorname{div}_g\operatorname{grad}_g=-d^*d
\]
on scalar functions. If \(D\) is flat, then \(d^2=(d^*)^2=0\).

We next fix the complex-geometric notation. Let \((M,g)\) be a Hermitian
manifold of complex dimension \(n\), let \((E,h)\) be a Hermitian holomorphic
vector bundle, and put
\[
F^{p,q}:=\Lambda^{p,q}T^*M\otimes E.
\]
If
\[
\nabla^E=\partial_E+\bar\partial_E,
\qquad
\Theta^{(E,h)}:=(\nabla^E)^2
\]
is the Chern connection and its curvature, we abbreviate
\(\bar\partial_E\) to \(\bar\partial\) and set
\[
\square_{p,q}
:=
\bar\partial\bar\partial^*
+\bar\partial^*\bar\partial.
\]
Let \(\nabla\) denote the induced connection on \(F^{p,q}\). Dependence of
a constant on \(E\) includes its metric, holomorphic structure, Chern
connection, and curvature.

Suppose that \(M\) is K\"ahler, with K\"ahler form \(\omega\). Let
\[
L_\omega\alpha:=\omega\wedge\alpha,\qquad
\Lambda:=L_\omega^*,
\qquad
\mathfrak{Ric}_{p,q}^E
:=
2\square_{p,q}-\nabla^*\nabla.
\]
The Akizuki--Nakano identity
\cite[Chapter~VII, Corollary~1.2]{D12} gives
\[
\mathfrak{Ric}_{p,q}^E
=
\mathfrak{Ric}_{p+q}^E+[i\Theta^{(E,h)},\Lambda],
\qquad
[A,B]:=A\circ B-B\circ A,
\]
where \(\mathfrak{Ric}_{p+q}^E\) is defined by
\eqref{equ:operator} for the induced Chern connection and restricted to
bidegree \((p,q)\). Consequently,
\[
\frac12\Delta|f|^2
=
\operatorname{Re}\langle
(-2\square_{p,q}+\mathfrak{Ric}_{p,q}^E)f,f\rangle
+|\nabla f|^2.
\]

Let \(L\) be a formally self-adjoint nonnegative Laplace-type operator on a
metric vector bundle \(F\) over a compact Riemannian manifold \(X\) with
smooth boundary. Its \emph{Dirichlet realization}, still denoted by \(L\),
has domain
\[
\operatorname{Dom}(L)=H^2(X,F)\cap H_0^1(X,F).
\]
This is the full-trace realization and is distinct from the relative,
absolute, and \(\bar\partial\)-Neumann realizations. If \(\ker L=\{0\}\),
then \(L^{-1}\) is the \emph{Dirichlet Green operator}; the kernels of
\(L^{-1}\) and \(e^{-tL}\) are called the \emph{Dirichlet Green kernel} and
\emph{Dirichlet heat kernel}, respectively.

Our kernel convention is
\[
K_T(x,y)\in\operatorname{Hom}(F_y,F_x)
\simeq F_x\otimes F_y^*,
\qquad
(Tf)(x)=\int_XK_T(x,y)f(y)\,dV(y).
\]
Unless stated otherwise, \(|K_T(x,y)|\) is the Hilbert--Schmidt norm;
fiberwise operator norms are denoted by \(\|\cdot\|_{\mathrm{op}}\). Thus,
if \(r_F=\operatorname{rank}F\),
\[
\|A\|_{\mathrm{op}}
\leq |A|
\leq\sqrt{r_F}\,\|A\|_{\mathrm{op}},
\qquad
|K_T(x,y)v|\leq |K_T(x,y)|\,|v|.
\]
Pairing with a kernel always means contraction in its dual factor.

For \(v\in F\), define its metric dual by
\[
v^*(w):=\langle w,v\rangle,
\qquad
Jv:=v^*.
\]
Thus \(J:F\to F^*\) is anti-linear and isometric; \(A^*\) denotes the
Hermitian adjoint of a fiber map \(A\). If
\[
L:C^\infty(U,F)\longrightarrow C^\infty(U,F')
\]
is a differential operator, its formal transpose
\[
L^\vee:C^\infty(U,(F')^*)\longrightarrow C^\infty(U,F^*)
\]
is characterized by
\[
\int_U(L^\vee\eta)(u)\,dV
=
\int_U\eta(Lu)\,dV
\]
for compactly supported smooth sections \(u\) and \(\eta\). If \(L^*\) is
the formal adjoint of \(L\), then
\[
L^\vee J_{F'}=J_FL^*.
\]

Finally, for an open set \(U\subset M^\circ\), using \(dV\) to identify
densities, we set
\[
\mathcal D(U,F):=C_c^\infty(U,F),
\qquad
\mathcal D'(U,F):=\bigl(\mathcal D(U,F^*)\bigr)',
\]
where the prime denotes the continuous complex-linear dual. For the trivial
complex line bundle, we write simply \(\mathcal D(U)\) and
\(\mathcal D'(U)\).

\section{Pointwise Estimates of Dirichlet Green Kernels}\label{sec:green form complex estimate}
This section proves Theorem~\ref{thm: green form complex boundary}.

Throughout this section, set
\[
F^{p,q}:=\Lambda^{p,q}T^*M\otimes E,
\qquad m:=\dim_{\mathbb R}M=2n.
\]
We use the same notation \(\square_{p,q}\) for its Dirichlet
realization on \(L^2(\Omega,F^{p,q})\), whose domain is
\[
\operatorname{Dom}(\square_{p,q})
=H^2(\Omega,F^{p,q})\cap H^1_0(\Omega,F^{p,q}).
\]
We write
\[
\delta(x):=d(x,\partial\Omega),\qquad x\in\overline\Omega.
\]
Constants in this section may depend on \((\overline\Omega,g,E)\), but not on
the points \(x,y\) or on the radius of a ball. Dependence on $E$ includes its
Hermitian and holomorphic structures and the associated Chern connection.
Unless otherwise indicated, norms of heat and Green kernels and their
bundle-valued derivatives are the pointwise Hermitian norms on the
corresponding tensor-product bundles. Fiberwise operator norms are written
explicitly as \(\|\cdot\|_{\mathrm{op}}\).

\subsection{The Dirichlet realization and its Green kernel}

\begin{lem}\label{lem:dirichlet-realization-kernel}
The kernel of the Dirichlet realization of \(\square_{p,q}\) is trivial. Consequently,
the first Dirichlet eigenvalue of \(\square_{p,q}\) is strictly positive.
\end{lem}
\begin{proof}
Let
\[
u\in H^2(\Omega,F^{p,q})\cap H^1_0(\Omega,F^{p,q}),
\qquad \square_{p,q}u=0.
\]
Since the full trace of \(u\) vanishes on \(\partial\Omega\), Green's
formula has no boundary contribution and gives
\[
0=(\square_{p,q}u,u)_{L^2(\Omega)}
=\|\bar\partial u\|_{L^2(\Omega)}^2
 +\|\bar\partial^*u\|_{L^2(\Omega)}^2.
\]
Hence
\[
\bar\partial u=0,\qquad \bar\partial^*u=0.
\]
Regard \(u\) as a section of the graded bundle
\[
\mathcal F^p:=\bigoplus_{s=0}^n\Lambda^{p,s}T^*M\otimes E
\]
and set
\[
\mathscr D_p:=\sqrt2\,(\bar\partial+\bar\partial^*).
\]
Then \(\mathscr D_pu=0\) in \(\Omega\), and \(\mathscr D_p\) is of Dirac
type. Extend the graded bundle and \(\mathscr D_p\) smoothly across
\(\partial\Omega\), and extend \(u\) by zero to the resulting neighborhood of
\(\overline\Omega\). Since
\(u\in H_0^1(\Omega,\mathcal F^p)\), its zero extension
\(\widetilde u\) belongs locally to \(H^1\).
Approximating \(u\) in \(H^1\) by compactly supported smooth sections shows
that
\[
\mathscr D_p\widetilde u=0
\]
in the distributional sense. On each component meeting \(\Omega\), the
extension vanishes on a nonempty exterior open set. Weak unique
continuation for Dirac-type operators \cite[Corollary~8.3]{BW93}, applied
componentwise, therefore gives \(\widetilde u=0\), and hence \(u=0\) in
\(\Omega\).

The Dirichlet realization is self-adjoint and has compact resolvent.
Therefore its spectrum is discrete, and triviality of its kernel implies
that its first eigenvalue is strictly positive.
\end{proof}

Define the Dirichlet Green operator by
\[
\mathcal G_{p,q}:=\square_{p,q}^{-1}.
\]
Let \(G_{p,q}(x,y)\) denote its Schwartz kernel. Thus
\[
G_{p,q}(x,y)\in
\operatorname{Hom}\bigl(F^{p,q}_y,F^{p,q}_x\bigr)
\simeq
F^{p,q}_x\otimes (F^{p,q}_y)^*,
\]
and
\[
(\mathcal G_{p,q}f)(x)
=
\int_\Omega G_{p,q}(x,y)f(y)\,dV(y),\qquad f\in L^2(\Omega,F^{p,q}).
\]
Let \(\square_{p,q}^{\vee}\) denote the formal transpose of
\(\square_{p,q}\) acting on the dual bundle \((F^{p,q})^*\), as defined in
Section~\ref{notation}. Whenever an operator carries the subscript \(y\), it
is understood to act on the \(y\)-variable, and hence on the corresponding
dual factor.

Let \(\{\phi_i\}_{i=1}^{\infty}\subset L^2(\Omega,F^{p,q})\)
be an orthonormal basis consisting of Dirichlet eigenforms of
\(\square_{p,q}\). Thus there exist eigenvalues
\[
0<\lambda_1\leq\lambda_2\leq\cdots,
\qquad
\lambda_i\to\infty,
\]
such that
\[
\begin{cases}
\square_{p,q}\phi_i=\lambda_i\phi_i
& \text{in }\Omega,\\
\phi_i|_{\partial\Omega}=0.
\end{cases}
\]
By elliptic boundary regularity,
\(\phi_i\in C^\infty(\overline\Omega,F^{p,q})\) for every $i$. The
Dirichlet heat kernel $H_{p,q}$ has the spectral representation
\[
H_{p,q}(t,x,y)=\sum_{i=1}^\infty
e^{-\lambda_it}\phi_i(x)\otimes\phi_i(y)^*,
\qquad t>0,\quad x,y\in\overline\Omega,
\]
with smooth convergence for $t$ bounded away from zero. As a Schwartz
kernel, the Green kernel has the distributional expansion
\[
G_{p,q}(x,y)
=
\sum_{i=1}^{\infty}
\frac{1}{\lambda_i}\,
\phi_i(x)\otimes\phi_i(y)^*.
\]
Here \(\phi_i(y)^*\in(F^{p,q}_y)^*\) is the metric dual of
\(\phi_i(y)\). Moreover, for $x\neq y$, the series converges pointwise,
\[
G_{p,q}(x,y)=\int_0^\infty H_{p,q}(t,x,y)\,dt,
\]
and the integral converges smoothly on compact subsets of
\((\overline\Omega\times\overline\Omega)\setminus\operatorname{diag}\).
By the spectral representation of $G_{p,q}$,
\[
\square_{p,q,x}G_{p,q}(x,y)
=
\delta_y(x)\,\operatorname{Id}_{F^{p,q}_y},\qquad
\square_{p,q,y}^{\vee}G_{p,q}(x,y)
=
\delta_x(y)\,\operatorname{Id}_{F^{p,q}_x}
\]
in the sense of distributions. Moreover, the full
Dirichlet boundary trace vanishes in each variable:
\[
G_{p,q}(x,\cdot)\big|_{\partial\Omega}=0,
\qquad
G_{p,q}(\cdot,y)\big|_{\partial\Omega}=0,
\]
in the trace sense.
Since the Dirichlet realization of \(\square_{p,q}\) is self-adjoint,
\[
G_{p,q}(x,y)=G_{p,q}(y,x)^*.
\]
Standard parabolic regularity gives
$H_{p,q}\in C^\infty((0,\infty)\times\overline\Omega\times\overline\Omega)$.
Elliptic interior and boundary regularity also imply that
\(G_{p,q}\in
C^\infty\!\left(
(\overline\Omega\times\overline\Omega)
\setminus\operatorname{diag}
\right);
\)
see, for example, \cite[Theorem~1.6.2]{S95}.

\subsection{A local first-order estimate}\

We prove the local first-order estimate used in the proof of
Theorem~\ref{thm: green form complex boundary} in this subsection.

For later use, set
\[
D_\Omega:=\operatorname{diam}(\overline\Omega),
\qquad
\kappa_0:=\max_{0\leq a,b\leq n}
\max\left\{0,-\min_{x\in\overline\Omega}
\lambda_{\min}\bigl(\mathfrak{Ric}_{a,b}^E(x)\bigr)\right\}.
\]
Thus $\mathfrak{Ric}_{a,b}^E\geq-\kappa_0 I\).

The next result extends a Euclidean estimate from \cite{Q25} to bundle-valued
forms on geodesic balls. Its proof adapts the Moser-iteration argument of Lu,
Qi S. Zhang, and Zhu \cite[Proposition~3.3]{LZZ21} to geodesic balls and
bundle-valued forms for the Dolbeault Laplacian.

\begin{thm}\label{thm:cheng yau general}
There exists a constant
\(C_{\mathrm{CY}}:=C_{\mathrm{CY}}(\overline\Omega,g,\kappa_0)>0\) such that,
for every geodesic ball \(B(x_0,r)\subset\Omega\) and every \(E\)-valued \((p,q)\)-form \(f\in C^0(\overline{B(x_0,r)},\Lambda^{p,q}T^*M\otimes E)\) satisfying
\(\square_{p,q}f=0\) in the distributional sense in \(B(x_0,r)\), one has
\begin{equation}\label{equ:local-first-order}
\sup_{B(x_0,r/2)}
\bigl(|\bar\partial f|^2+|\bar\partial^*f|^2\bigr)
\leq \frac{C_{\mathrm{CY}}}{r^2}
\sup_{B(x_0,r)}|f|^2.
\end{equation}
\end{thm}
\begin{proof}
We follow \cite[Proposition~3.3]{LZZ21}. Elliptic regularity implies that $f$
is smooth in $B(x_0,r)$.
Fix a relatively compact neighborhood \(U\) of \(\overline\Omega\).
After extending the metric outside \(U\), we may regard \(U\) as a subset
of a complete Riemannian manifold without changing the metric on \(U\).
Choose
\[
0<r_0\leq \min\{1,D_\Omega\}
\]
so small that \(B(z,4r_0)\Subset U\) for every
\(z\in\overline\Omega\), and that normal coordinates are uniformly
controlled on these balls. Uniform equivalence of the metric coefficients
and the Euclidean metric in these coordinates, and \cite[Theorem~10.4, p.~445]{S92} give a constant
\(C_0:=C_{0}(\overline\Omega,g)\geq1\) such that the following estimates hold for every \(z\in\overline\Omega\),
\(0<\rho\leq r_0\), and \(\psi\in C_c^\infty(B(z,\rho))\):
\begin{align}
 C_0^{-1}\rho^m
 &\leq \operatorname{Vol}B(z,\rho)\leq C_0\rho^m,
 \label{equ:local-volume-single-constant}\\
 \left(\fint_{B(z,\rho)}|\psi|^{2\chi}\,dV\right)^{1/\chi}
 &\leq C_0\rho^2
 \fint_{B(z,\rho)}|\nabla\psi|^2\,dV,
 \label{equ:local-sobolev-single-constant}
\end{align}
where
\[
\chi:=\begin{cases}
\dfrac{m}{m-2},&m>2,\\[1mm]
2,&m=2.
\end{cases}
\]

We first assume \(0<r\leq r_0\). Set
\[
u_+:=\bar\partial f,
\qquad
u_-:=\bar\partial^*f,
\qquad
v:=|u_+|^2+|u_-|^2,
\]
where a nonexistent term is omitted when \(q=n\) or \(q=0\). Since
\(\bar\partial^2=(\bar\partial^*)^2=0\), one has on smooth forms
\begin{align*}
\square_{p,q+1}\bar\partial
 =\bar\partial\square_{p,q},\qquad
\square_{p,q-1}\bar\partial^*
 =\bar\partial^*\square_{p,q}.
\end{align*}
Consequently,
\[
\square_{p,q+1}u_+=0,
\qquad
\square_{p,q-1}u_-=0.
\]
We begin with an explicit local mean-value estimate for either
section \(u=u_+\) or \(u=u_-\). The Bochner--Weitzenb\"ock formula gives
\begin{equation}\label{equ:bochner-u-explicit}
\Delta|u|^2
\geq -2\kappa_0|u|^2+2|\nabla u|^2.
\end{equation}
For \(\varepsilon>0\), set
\(w_\varepsilon=(|u|^2+\varepsilon^2)^{1/2}\). The chain rule for the
Laplace--Beltrami operator gives
\begin{align*}
\Delta w_\varepsilon
&=\frac{\Delta|u|^2}{2w_\varepsilon}
 -\frac{\bigl|d(|u|^2)\bigr|^2}{4w_\varepsilon^3}\geq -\kappa_0\frac{|u|^2}{w_\varepsilon}
 +\frac{|\nabla u|^2}{w_\varepsilon}
 -\frac{\bigl|d(|u|^2)\bigr|^2}{4w_\varepsilon^3}.
\end{align*}
However, Kato's inequality yields
\[
\bigl|d(|u|^2)\bigr|^2
\leq4|u|^2|\nabla u|^2.
\]
Therefore
\[
\Delta w_\varepsilon
\geq-\kappa_0w_\varepsilon
+\frac{\varepsilon^2}{w_\varepsilon^3}|\nabla u|^2
\geq-\kappa_0w_\varepsilon.
\]
We first prove that every positive continuous weak supersolution
\(w\) of \(\Delta w\geq-\kappa_0w\) satisfies
\begin{equation}\label{equ:explicit-local-mean-value}
\sup_{B(x_0,r/2)}w^2
\leq C_{\mathrm M}r^{-m}
\int_{B(x_0,3r/4)}w^2\,dV.
\end{equation}

Let \(p\geq2\), let
\(0\leq\eta\in C_c^\infty(B(x_0,r))\), and test
\(\Delta w\geq-\kappa_0w\) against \(\eta^2w^{p-1}\). Integration
by parts gives
\[
(p-1)\int_{B(x_0,r)}\eta^2w^{p-2}|\nabla w|^2
\leq2\int_{B(x_0,r)}\eta w^{p-1}|\nabla w|\,|\nabla\eta|
+\kappa_0\int_{B(x_0,r)}\eta^2w^p.
\]
Put \(\varphi:=w^{p/2}\). Since
\[
|\nabla \varphi|=\frac p2w^{p/2-1}|\nabla w|,
\]
we obtain
\begin{equation}\label{equ:moser-energy-first}
\frac{4(p-1)}{p^2}
\int_{B(x_0,r)}\eta^2|\nabla \varphi|^2
\leq \frac4p\int_{B(x_0,r)}\eta \varphi|\nabla \varphi||\nabla\eta|
 +\kappa_0\int_{B(x_0,r)}\eta^2\varphi^2.
\end{equation}
Set
\[
I:=\int_{B(x_0,r)}\eta^2|\nabla \varphi|^2,\qquad
J:=\int_{B(x_0,r)}|\nabla\eta|^2\varphi^2,\qquad
K:=\int_{B(x_0,r)}\eta^2\varphi^2.
\]
By Young's inequality,
\[
\frac4p\sqrt{IJ}
\leq\frac{2(p-1)}{p^2}I+\frac{2}{p-1}J.
\]
Substitution into \eqref{equ:moser-energy-first} yields
\[
I\leq\frac{p^2}{(p-1)^2}J
+\frac{p^2}{2(p-1)}\kappa_0K.
\]
Consequently, for any $p\geq 2$,
\begin{align}\label{equ:energy}
\int_{B(x_0,r)}|\nabla(\eta w^{p/2})|^2
&\leq2I+2J\leq
\left(\frac{2p^2}{(p-1)^2}+2\right)J
+\frac{p^2}{p-1}\kappa_0K\nonumber\\
&\leq4p^2\int_{B(x_0,r)}
\bigl(|\nabla\eta|^2+\kappa_0\eta^2\bigr)w^p.
\end{align}

Set
\[
\rho_j:=\frac r2+\frac{r}{2^{j+2}},
\qquad
p_j:=2\chi^j,
\qquad j=0,1,\ldots,
\]
and choose \(0\leq \eta_j\in C_c^\infty(B(x_0,\rho_j))\) equal to one on
\(B(x_0,\rho_{j+1})\), with
\[
|\nabla\eta_j|\leq \frac{2^{j+4}}r.
\]
Define
\[
Y_j:=\left(
\fint_{B(x_0,\rho_j)}w^{p_j}\,dV
\right)^{1/p_j}.
\]
Using \eqref{equ:local-sobolev-single-constant} and \eqref{equ:energy}
we obtain
\begin{align*}
Y_{j+1}^{p_j}
&\leq
4C_0
\left(
\frac{\operatorname{Vol}B(x_0,\rho_j)}
     {\operatorname{Vol}B(x_0,\rho_{j+1})}
\right)^{1/\chi}
p_j^2\rho_j^2
\left(\frac{2^{2j+8}}{r^2}+\kappa_0\right)
Y_j^{p_j}.
\end{align*}
Taking the \(p_j\)-th root and then using the volume-ratio bound
\[
\frac{\operatorname{Vol}B(x_0,\rho_j)}
     {\operatorname{Vol}B(x_0,\rho_{j+1})}
\leq \left(\frac32\right)^m C_0^2,
\]
we obtain
\begin{align}
Y_{j+1}\leq
\left[
576C_0\left(\left(\frac32\right)^mC_0^2\right)^{1/\chi}
(1+\kappa_0r_0^2)p_j^2 4^j
\right]^{1/p_j}Y_j.
\label{equ:explicit-moser-step}
\end{align}
A direct calculation gives
\[
\sum_{j=0}^\infty\frac1{p_j}
 =\frac{\chi}{2(\chi-1)},
\qquad
\sum_{j=0}^\infty\frac{j}{p_j}
 =\frac{\chi}{2(\chi-1)^2}.
\]
Set
\[
\mathscr A:=
576C_0\left(\left(\frac32\right)^mC_0^2\right)^{1/\chi}
(1+\kappa_0r_0^2).
\]
Since \(p_j=2\chi^j\), iteration of
\eqref{equ:explicit-moser-step} gives
\begin{align*}
\limsup_{j\to\infty}Y_j
&\leq
\mathscr A^{\sum_{j=0}^\infty 1/p_j}
\prod_{j=0}^\infty p_j^{2/p_j}4^{j/p_j}Y_0\\
&=
\mathscr A^{\frac{\chi}{2(\chi-1)}}
2^{\frac{\chi}{\chi-1}+
  \frac{\chi}{(\chi-1)^2}}
\chi^{\frac{\chi}{(\chi-1)^2}}Y_0.
\end{align*}
Since \(\rho_j\to r/2\) as $j\rightarrow\infty$,
\(B(x_0,r/2)\subset B(x_0,\rho_j)\), and \(w\) is continuous, the
standard \(L^p\)-limit on these nested balls gives
\[
\sup_{B(x_0,r/2)}w
\leq\liminf_{j\to\infty}Y_j.
\]
Furthermore, the lower volume bound gives
\[
Y_0^2
=\fint_{B(x_0,3r/4)}w^2\,dV
\leq\left(\frac43\right)^mC_0r^{-m}
\int_{B(x_0,3r/4)}w^2\,dV.
\]
Consequently, \eqref{equ:explicit-local-mean-value} holds with
\begin{equation}\label{equ:explicit-local-mean-value-constant}
C_{\mathrm M}:=
\left(\frac43\right)^m C_0
2^{\frac{2\chi}{\chi-1}+\frac{2\chi}{(\chi-1)^2}}
\chi^{\frac{2\chi}{(\chi-1)^2}}
\mathscr A^{\frac{\chi}{\chi-1}}.
\end{equation}
Applying \eqref{equ:explicit-local-mean-value} separately to the regularized
norms of $u_+$ and $u_-$, adding the resulting inequalities, and then letting
\(\varepsilon\downarrow0\), gives
\begin{equation}\label{equ:explicit-moser-v}
\sup_{B(x_0,r/2)}v
\leq C_{\mathrm M}r^{-m}
\int_{B(x_0,3r/4)}v\,dV.
\end{equation}

It remains to estimate the integral of \(v\). Choose a smooth radial
cut-off \(0\leq \eta\in C_c^\infty(B(x_0,r))\), equal to one on
\(B(x_0,3r/4)\), and satisfying \(|\nabla\eta|\leq8/r\). Then the assumption \(\square_{p,q}f=0\) gives
\begin{align*}
0={}&\int_{B(x_0,r)}\langle \square_{p,q}f,\eta^2f\rangle dV\\
={}&\int_{B(x_0,r)}\eta^2
\bigl(|\bar\partial f|^2+|\bar\partial^*f|^2\bigr)\,dV+2\operatorname{Re}\int_{B(x_0,r)}\eta
\langle\bar\partial f,\bar\partial\eta\wedge f\rangle\,dV\\
{}&-2\operatorname{Re}\int_{B(x_0,r)}\eta
\left\langle \bar\partial\eta\wedge \bar\partial^*f,
f\right\rangle\,dV.
\end{align*}
By the Cauchy--Schwarz inequality,
\begin{align*}
\int_{B(x_0,r)}\eta^2v\,dV&\leq 2\int_{B(x_0,r)}\eta|\nabla\eta||f|(|\bar\partial f|+|\bar\partial^*f|)dV\\
&\leq2\sqrt2\int_{B(x_0,r)}\eta|\nabla\eta||f|v^{1/2}\,dV\\
&\leq\frac12\int_{B(x_0,r)}\eta^2v\,dV
 +4\int_{B(x_0,r)}|\nabla\eta|^2|f|^2\,dV.
\end{align*}
Consequently,
\begin{equation}\label{equ:explicit-caccioppoli}
\int_{B(x_0,3r/4)}v\,dV
\leq\frac{512}{r^2}
\int_{B(x_0,r)}|f|^2\,dV.
\end{equation}
Combining \eqref{equ:explicit-moser-v},
\eqref{equ:explicit-caccioppoli}, and the upper volume bound in
\eqref{equ:local-volume-single-constant}, we obtain
\[
\sup_{B(x_0,r/2)}v
\leq\frac{512C_0C_{\mathrm M}}{r^2}
\sup_{B(x_0,r)}|f|^2
\qquad(0<r\leq r_0).
\]

Finally, suppose \(r>r_0\). For every \(z\in B(x_0,r/2)\), the ball
\(B(z,r_0/2)\) is contained in \(B(x_0,r)\). Applying the already proved
small-radius estimate on this ball gives
\[
v(z)\leq\frac{4\cdot512C_0C_{\mathrm M}}{r_0^2}
\sup_{B(x_0,r)}|f|^2\leq
\frac{512C_0C_{\mathrm M}}{r^2}
\frac{4D_\Omega^2}{r_0^2}
\sup_{B(x_0,r)}|f|^2.
\]
Here we used $r\leq D_\Omega$.
Thus \eqref{equ:local-first-order} holds with
\begin{equation}\label{equ:explicit-CY-constant}
C_{\mathrm{CY}}
:=512C_0C_{\mathrm M}
\max\left\{1,\frac{4D_\Omega^2}{r_0^2}\right\}.
\end{equation}
\end{proof}

\begin{cor}\label{thm:cheng yau cor}
Under the hypotheses of Theorem~\ref{thm:cheng yau general},
\[
\max\left\{
\sup_{B(x_0,r/2)}|\bar\partial f|,
\sup_{B(x_0,r/2)}|\bar\partial^*f|
\right\}
\leq \frac{\sqrt{C_{\mathrm{CY}}}}{r}
\sup_{B(x_0,r)}|f|.
\]
The same estimate, with the same constant, holds for the formal transpose
on the dual bundle and hence for the \(y\)-sections of the Green kernel away
from the pole.
\end{cor}
\begin{proof}
The first assertion follows by taking square roots in
\eqref{equ:local-first-order}. Let \(\square_{p,q}^{\vee}\) denote the
formal transpose acting on \(F^{p,q*}\). Its principal metric and its local
Sobolev and volume constants are unchanged, while its Weitzenb\"ock
endomorphism is the fiberwise dual transpose of the original one. It therefore
has the same pointwise eigenvalues, so the same \(\kappa_0\) and
\(C_{\mathrm{CY}}\) apply, and the last claim follows.
\end{proof}

\subsection{Long-time Dirichlet heat kernel estimates}\ 

For the remaining estimates in this section, let \(H_0\) be the scalar
Dirichlet heat kernel associated with \(e^{t\Delta}\). Thus,
\[
(\partial_t-\Delta_x)H_0(t,x,y)=0
\]
with zero boundary values.

\begin{lem}[Scalar Dirichlet heat kernel]\label{lem:scalar-dirichlet-heat}
There is a constant $A_H:=A_H(\overline\Omega,g)\geq1$ such that, for
$0<t\leq1$ and $x,y\in\Omega$,
\begin{equation}\label{equ:scalar-full-boundary-gaussian}
H_0(t,x,y)\leq A_Ht^{-m/2}
\min\left\{1,\frac{\delta(x)}{\sqrt t}\right\}
\min\left\{1,\frac{\delta(y)}{\sqrt t}\right\}
\exp\left(-\frac{d(x,y)^2}{A_Ht}\right).
\end{equation}
\end{lem}
This is the standard short-time boundary estimate for the scalar Dirichlet
heat kernel; see \cite{D87,D89,Z02}.
We repeatedly use the following elementary consequence of a heat-trace bound.
Suppose that
\[
0<\nu_1\leq \nu_2\leq\cdots
\]
and
\[
\sum_{j=1}^{\infty} e^{-t\nu_j}\leq A t^{-m/2},
\qquad 0<t\leq 1.
\]
Then, for every \(j\geq1\),
\begin{equation}\label{equ:explicit-heat-trace-eigenvalue}
\nu_j\geq \frac{m}{2e}A^{-2/m}j^{2/m}-\frac{m}{2}.
\end{equation}
Indeed, by the monotonicity of \((\nu_j)\),
\[
j e^{-t\nu_j}
\leq \sum_{\ell=1}^{j}e^{-t\nu_\ell}
\leq A t^{-m/2},
\qquad 0<t\leq1.
\]
If \(\nu_j\geq m/2\), choosing
\[
t=\frac{m}{2\nu_j}\leq1
\]
gives
\[
\nu_j\geq \frac{m}{2e}A^{-2/m}j^{2/m}.
\]
In particular, \eqref{equ:explicit-heat-trace-eigenvalue} follows.
If instead \(\nu_j<m/2\), taking \(t=1\) yields
\[
j e^{-\nu_j}\leq A,
\]
so that
\[
j\leq A e^{\nu_j}<A e^{m/2}.
\]
Consequently,
\[
\frac{m}{2e}A^{-2/m}j^{2/m}-\frac{m}{2}<0<\nu_j,
\]
which again proves \eqref{equ:explicit-heat-trace-eigenvalue}.

We also use the following elementary estimate. For \(a>0\), set
\begin{equation*}
M_a:=\sup_{s\geq0}e^{-s/4}(1+s)^a\leq e(4a)^a.
\end{equation*}
It follows that
\begin{align}
\sum_{j=1}^{\infty}e^{-\nu_j/2}(1+\nu_j)^a
&\leq M_a\sum_{j=1}^{\infty}e^{-\nu_j/4}\leq 2^{m+2}(4a)^a A.
\label{equ:spectral-sum-bound}
\end{align}

\begin{lem}[Full-trace Dirichlet semigroup domination]
\label{lem:full-dirichlet-semigroup-domination}
For every $t>0$ and $s\in L^2(\Omega,F^{p,q})$,
\begin{equation}\label{equ:dirichlet-semigroup-domination}
\bigl|e^{-t\square_{p,q}}s\bigr|
\leq e^{\kappa_0t/2}e^{(t/2)\Delta}|s|
\qquad\text{a.e. on }\Omega.
\end{equation}
Consequently, the heat kernels satisfy
\begin{equation}\label{equ:dirichlet-kernel-op-domination}
\|H_{p,q}(t,x,y)\|_{\mathrm{op}}
\leq e^{\kappa_0t/2}H_0(t/2,x,y).
\end{equation}
\end{lem}
Both assertions are the Dirichlet semigroup-domination estimate of
Hess--Schrader--Uhlenbrock; see \cite{HSU77} and
\cite[Theorem~4.3]{DL82}.

\begin{lem}\label{lem:bundle-eigenfunction}
Let \(\phi\) be an \(L^2\)-normalized Dirichlet eigenform of
\(\square_{p,q}\):
\[
\square_{p,q}\phi=\lambda\phi,
\qquad \phi|_{\partial\Omega}=0.
\]
Then
\begin{equation}\label{equ:bundle-eigenfunction}
\|\phi\|_{L^\infty(\Omega)}^2
\leq 2^{m/2+2}e^{\kappa_0/2}A_H(1+\lambda)^{m/2},
\end{equation}
and
\begin{equation}\label{equ:bundle-eigenfunction-boundary}
|\phi(x)|^2
\leq 2^{m/2+3}e^{\kappa_0/2}A_H\delta(x)^2
(1+\lambda)^{m/2+1},
\qquad x\in\Omega.
\end{equation}
\end{lem}
\begin{proof}
By Lemma~\ref{lem:full-dirichlet-semigroup-domination}, for any $s>0$,
\[
\|H_{p,q}(s,x,y)\|_{\mathrm{op}}
\leq e^{\kappa_0s/2}H_0(s/2,x,y).
\]
The spectral definition of $H_{p,q}$ gives
\[
e^{-t\lambda}\phi(x)
=\int_\Omega H_{p,q}(t,x,z)\phi(z)\,dV(z).
\]
Taking norms and applying semigroup domination gives the first inequality
below; the scalar Cauchy--Schwarz inequality and the semigroup identity give
the second:
\begin{align*}
e^{-2t\lambda}|\phi(x)|^2
&\leq e^{\kappa_0t}
\left(\int_\Omega H_0(t/2,x,z)|\phi(z)|\,dV(z)\right)^2\leq e^{\kappa_0t}H_0(t,x,x).
\end{align*}
Set \(t=[2(1+\lambda)]^{-1}\). The scalar estimate
\eqref{equ:scalar-full-boundary-gaussian}, first without and then with its two
diagonal boundary factors, gives \eqref{equ:bundle-eigenfunction} and
\eqref{equ:bundle-eigenfunction-boundary}.
\end{proof}

\begin{lem}\label{lem:bundle-spectrum}
Let \(0<\lambda_1\leq\lambda_2\leq\cdots\) be the Dirichlet spectrum
of \(\square_{p,q}\), and set
\begin{equation}\label{equ:explicit-bundle-heat-trace-constant}
A_{p,q}
:=
2^{3m/2}\operatorname{rank}(E)e^{\kappa_0/2}
A_H\operatorname{Vol}(\Omega).
\end{equation}
Then
\begin{equation}\label{equ:bundle-DL-lower}
\lambda_j\geq \frac{m}{2e}A_{p,q}^{-2/m}j^{2/m}-\frac m2.
\end{equation}
Let
\[
\lambda_0:=\min_{0\leq a,b\leq n}
\lambda_1(\square_{a,b})>0.
\]
For \(t\geq1\),
\begin{equation}
|H_{p,q}(t,x,y)|
\leq C_{\mathrm{long}}e^{-\lambda_0t/2},\ 
|H_{p,q}(t,x,y)|\leq C_{\mathrm{long}}\delta(y)e^{-\lambda_0t/2},
\label{equ:bundle-long-boundary}
\end{equation}
\begin{equation}
|H_{p,q}(t,x,y)|
\leq C_{\mathrm{long}}\delta(x)\delta(y)
e^{-\lambda_0t/2},
\label{equ:bundle-long-two-boundary}
\end{equation}
where
\begin{align}
C_{\mathrm{long}}:=2^{2m+6}(m+2)^{(m+2)/2}e^{\kappa_0/2}A_H A_{p,q}.
\label{equ:explicit-bundle-two-boundary-long-constant}
\end{align}
\end{lem}
\begin{proof}
For \(0<t\leq1\), semigroup domination
\cite[Theorem~4.3]{DL82} on the diagonal gives
\begin{align*}
\sum_{j=1}^\infty e^{-t\lambda_j}
&=\int_\Omega\operatorname{tr}H_{p,q}(t,x,x)\,dV(x)\\
&\leq\operatorname{rank}(F^{p,q})e^{\kappa_0t/2}
\int_\Omega H_0(t/2,x,x)\,dV(x)\leq A_{p,q}t^{-m/2}.
\end{align*}
Here
\(\operatorname{rank}(F^{p,q})
=\binom np\binom nq\operatorname{rank}(E)
\leq2^m\operatorname{rank}(E)\), which explains the choice of
$A_{p,q}$.
Thus \eqref{equ:bundle-DL-lower} follows from
\eqref{equ:explicit-heat-trace-eigenvalue}. Lemma
\ref{lem:dirichlet-realization-kernel} implies
\(\lambda_0>0\).

For \(t\geq1\), Lemma~\ref{lem:bundle-eigenfunction}, the spectral expansion of $H_{p,q}$
and \eqref{equ:spectral-sum-bound} give
\begin{align*}
|H_{p,q}(t,x,y)|
&\leq 2^{m/2+2}e^{\kappa_0/2}A_He^{-\lambda_0t/2}
\sum_{j=1}^\infty
e^{-\lambda_j/2}(1+\lambda_j)^{m/2}\leq C_{\mathrm{long}}e^{-\lambda_0t/2}.
\end{align*}
\begin{align*}
|H_{p,q}(t,x,y)|
&\leq 2^{m/2+3}e^{\kappa_0/2}A_He^{-\lambda_0t/2}
\sum_{j=1}^\infty
e^{-\lambda_j/2}(1+\lambda_j)^{(m+1)/2}\delta(y)\\
&\leq C_{\mathrm{long}}\delta(y)e^{-\lambda_0t/2}.
\end{align*}
\begin{align*}
|H_{p,q}(t,x,y)|
&\leq 2^{m/2+3}e^{\kappa_0/2}A_H
\delta(x)\delta(y)e^{-\lambda_0t/2}
\sum_{j=1}^\infty e^{-\lambda_j/2}
(1+\lambda_j)^{(m+2)/2}\\
&\leq C_{\mathrm{long}}\delta(x)\delta(y)
e^{-\lambda_0t/2}.
\end{align*}
This proves the lemma.
\end{proof}

\subsection{Integration of the heat kernel estimates}\ 

 In this subsection, we first integrate the zero-order heat kernel bounds and then estimate the
first-order kernels, completing the proof of
Theorem~\ref{thm: green form complex boundary}.
\begin{lem}
\label{lem:bundle-green-from-heat}
For distinct \(x,y\in\Omega\), put \(r=d(x,y)\).  If \(m>2\), then
\[
|G_{p,q}(x,y)|\leq C_{\mathrm G}^{(m)}r^{2-m},
\]
whereas, if \(m=2\), then
\[
|G_{p,q}(x,y)|\leq C_{\mathrm G}^{(2)}(1+|\ln r|).
\]
For every \(m\geq2\), one also has
\begin{align*}
|G_{p,q}(x,y)|\leq C_{\partial\mathrm G}\delta(y)r^{1-m},\ 
|G_{p,q}(x,y)|\leq C_{\partial^2\mathrm G}\delta(x)\delta(y)r^{-m},
\end{align*}
and the constants are given explicitly in
\eqref{equ:explicit-green-constant-high-dim},
\eqref{equ:explicit-green-constant-dim2}, 
\eqref{equ:explicit-green-boundary-constant}, and 
\eqref{equ:explicit-green-two-boundary-constant}.
\end{lem}
\begin{proof}
By the definition of \(\kappa_0\),
\[
\mathfrak{Ric}_{p,q}^E\geq-\kappa_0I.
\]
Lemma~\ref{lem:full-dirichlet-semigroup-domination} and the norm convention
in Section~\ref{notation} give
\[
|H_{p,q}(t,x,y)|
\leq \sqrt{\operatorname{rank}F^{p,q}}\,
e^{\kappa_0t/2}H_0(t/2,x,y).
\]
Consequently, for $0<t\leq 1$, applying \eqref{equ:scalar-full-boundary-gaussian} at time \(t/2\),
we obtain
\begin{equation}
|H_{p,q}(t,x,y)|
\leq C_{\mathrm{sh}}\min\left\{1,\frac{\delta(y)}{\sqrt{t}}\right\}t^{-m/2}
\exp\left(-\frac{d(x,y)^2}{A_Ht}\right),
\label{equ:bundle-short-boundary}
\end{equation}
\begin{align}
|H_{p,q}(t,x,y)|
&\leq C_{\mathrm{sh}}
\min\left\{1,\frac{\delta(x)}{\sqrt{t}}\right\}
\min\left\{1,\frac{\delta(y)}{\sqrt{t}}\right\}t^{-m/2}
\exp\left(-\frac{d(x,y)^2}{A_Ht}\right).
\label{equ:bundle-short-two-boundary}
\end{align}
Here
\begin{align}
C_{\mathrm{sh}}:=
2^{m/2+1}\sqrt{\operatorname{rank}F^{p,q}}\,
e^{\kappa_0/2}A_H.
\label{equ:explicit-bundle-short-two-boundary-constant}
\end{align}

Set
\begin{align*}
L_{\mathrm{long}}:=\frac{2C_{\mathrm{long}}}{\lambda_0}e^{-\lambda_0/2}.
\end{align*}
By Lemma \ref{lem:bundle-spectrum}, we have
\begin{align}
\int_1^\infty|H_{p,q}(t,x,y)|\,dt\leq L_{\mathrm{long}}\cdot
\min\left\{1,\delta(x),\delta(y),\delta(x)\delta(y)\right\}.
\label{equ:integrated-long-two-boundary}
\end{align}

It is clear that 
\[
\int_0^1 t^{-m/2}
\exp\left(-\frac{r^2}{A_Ht}\right)\,dt
\leq
\begin{cases}
A_H^{m/2-1}
\Gamma\left(\dfrac m2-1\right)r^{2-m},
& m>2,\\[2mm]
1+\ln A_H+2|\ln r|,
& m=2.
\end{cases}
\]
If $m>2$, then
\begin{align*}
  |G_{p,q}(x,y)|&\leq\int_0^{\infty}|H_{p,q}(t,x,y)|dt\leq C_{\mathrm{sh}}A_H^{m/2-1}
\Gamma\left(\frac m2-1\right)r^{2-m}+L_{\mathrm{long}}\\
&\leq C_{\mathrm G}^{(m)}r^{2-m},
\end{align*}
where
\begin{equation}\label{equ:explicit-green-constant-high-dim}
C_{\mathrm G}^{(m)}
:=C_{\mathrm{sh}}A_H^{m/2-1}
\Gamma\left(\frac m2-1\right)
+L_{\mathrm{long}}D_\Omega^{m-2}.
\end{equation}
If $m=2$, then
\begin{align*}
  |G_{p,q}(x,y)|&\leq\int_0^{\infty}|H_{p,q}(t,x,y)|dt\leq C_{\mathrm{sh}}(1+\ln A_H+2|\ln r|)+L_{\mathrm{long}}\\
  &\leq C_{\mathrm G}^{(2)}(1+|\ln r|),
\end{align*}
where
\begin{equation}\label{equ:explicit-green-constant-dim2}
C_{\mathrm G}^{(2)}
:=C_{\mathrm{sh}}
\bigl(3+\ln A_H\bigr)+L_{\mathrm{long}}.
\end{equation}

Note that
\begin{align*}
\int_0^1t^{-(m+1)/2}
e^{-r^2/(A_Ht)}\,dt\leq A_H^{(m-1)/2}
\Gamma\left(\frac{m-1}{2}\right)r^{1-m}.
\end{align*}
The same argument gives
\begin{equation}\label{equ:bundle-green-boundary}
|G_{p,q}(x,y)|
\leq C_{\partial\mathrm G}\delta(y)r^{1-m},
\end{equation}
where
\begin{equation}\label{equ:explicit-green-boundary-constant}
C_{\partial\mathrm G}
:=C_{\mathrm{sh}}A_H^{(m-1)/2}
\Gamma\left(\frac{m-1}{2}\right)
+L_{\mathrm{long}}D_\Omega^{m-1}.
\end{equation}
Finally,
\[
\int_0^1t^{-(m+2)/2}e^{-r^2/(A_Ht)}\,dt
\leq A_H^{m/2}\Gamma\left(\frac m2\right)r^{-m}.
\]
Keeping both boundary factors in
\eqref{equ:bundle-short-two-boundary} and using
\eqref{equ:integrated-long-two-boundary}, we obtain
\begin{equation}\label{equ:bundle-green-two-boundary}
|G_{p,q}(x,y)|
\leq C_{\partial^2\mathrm G}\delta(x)\delta(y)r^{-m},
\end{equation}
where
\begin{equation}\label{equ:explicit-green-two-boundary-constant}
C_{\partial^2\mathrm G}
:=C_{\mathrm{sh}}A_H^{m/2}\Gamma\left(\frac m2\right)
+L_{\mathrm{long}}D_\Omega^m.
\end{equation}
\end{proof}
\begin{proof}[Proof of Theorem~\ref{thm: green form complex boundary}]
 Fix distinct
\(x,y\in\overline\Omega\) and set \(r=d(x,y)\). Since \(G_{p,q}\) is smooth
on\(
(\overline\Omega\times\overline\Omega)\setminus\operatorname{diag},\)
we may assume $x,y\in\Omega$. Lemma~\ref{lem:bundle-green-from-heat}
establishes the zero-order estimates, including the estimate with both
boundary factors. It also gives the estimate with the factor
\(\delta(y)\). The Hermitian
symmetry \(G_{p,q}(x,y)=G_{p,q}(y,x)^*\) gives the same estimate with
\(\delta(x)\); taking the smaller of the two proves the first estimate in
(ii), while \eqref{equ:bundle-green-two-boundary} proves the second. It remains to
estimate the first-order terms.

Suppose first that \(\delta(y)\leq r\).
By Corollary~\ref{thm:cheng yau cor},
\[
\max\bigl\{|\bar\partial_yG_{p,q}(x,y)|,
|\bar\partial_y^*G_{p,q}(x,y)|\bigr\}
\leq\frac{2\sqrt{C_{\mathrm{CY}}}}{\delta(y)}
\sup_{z\in B(y,\delta(y)/2)}|G_{p,q}(x,z)|.
\]
For any \(z\in B(y,\delta(y)/2)\), the triangle inequality gives
\[
d(x,z)\geq r-\frac{\delta(y)}2\geq\frac r2,
\qquad
\delta(z)\leq\delta(y)+d(y,z)\leq\frac32\delta(y).
\]
Consequently, Lemma~\ref{lem:bundle-green-from-heat} gives
\[
\sup_{z\in B(y,\delta(y)/2)}|G_{p,q}(x,z)|
\leq\frac32\,2^{m-1}C_{\partial\mathrm G}
\delta(y)r^{1-m}.
\]
Hence
\begin{equation}\label{equ:first-derivative-case-one-explicit}
\max\bigl\{|\bar\partial_yG_{p,q}(x,y)|,
|\bar\partial_y^*G_{p,q}(x,y)|\bigr\}
\leq3\cdot2^{m-1}\sqrt{C_{\mathrm{CY}}}
C_{\partial\mathrm G}r^{1-m}.
\end{equation}
Keeping both boundary factors instead gives
\[
\sup_{z\in B(y,\delta(y)/2)}|G_{p,q}(x,z)|
\leq \frac32\,2^m C_{\partial^2\mathrm G}
\delta(x)\delta(y)r^{-m}.
\]
Consequently,
\begin{equation}\label{equ:first-derivative-weighted-case-one}
\max\bigl\{|\bar\partial_yG_{p,q}(x,y)|,
|\bar\partial_y^*G_{p,q}(x,y)|\bigr\}
\leq3\cdot2^m\sqrt{C_{\mathrm{CY}}}
C_{\partial^2\mathrm G}\delta(x)r^{-m}.
\end{equation}

Suppose next that \(\delta(y)>r\). Then \(B(y,r/2)\Subset\Omega\) and
\(x\notin B(y,r/2)\). Then Corollary~\ref{thm:cheng yau cor} gives
\begin{equation}\label{equ:1 order}
\max\bigl\{|\bar\partial_yG_{p,q}(x,y)|,
|\bar\partial_y^*G_{p,q}(x,y)|\bigr\}
\leq\frac{2\sqrt{C_{\mathrm{CY}}}}r
\sup_{z\in B(y,r/2)}|G_{p,q}(x,z)|.
\end{equation}
For \(z\in B(y,r/2)\),
\[
\frac r2\leq d(x,z)\leq\frac{3r}{2}.
\]
By the Hermitian symmetry and
\eqref{equ:bundle-green-boundary},
\[
\sup_{z\in B(y,r/2)}|G_{p,q}(x,z)|
\leq2^{m-1}C_{\partial\mathrm G}\delta(x)r^{1-m}.
\]
It follows from \eqref{equ:1 order} that
\begin{equation}\label{equ:first-derivative-weighted-case-two}
\max\bigl\{|\bar\partial_yG_{p,q}(x,y)|,
|\bar\partial_y^*G_{p,q}(x,y)|\bigr\}
\leq2^m\sqrt{C_{\mathrm{CY}}}
C_{\partial\mathrm G}\delta(x)r^{-m}.
\end{equation}
If \(m>2\), then Lemma~\ref{lem:bundle-green-from-heat} implies
\[
\sup_{z\in B(y,r/2)}|G_{p,q}(x,z)|
\leq2^{m-2}C_{\mathrm G}^{(m)}r^{2-m},
\]
and therefore
\begin{equation}\label{equ:first-derivative-case-two-high}
\max\bigl\{|\bar\partial_yG_{p,q}(x,y)|,
|\bar\partial_y^*G_{p,q}(x,y)|\bigr\}
\leq2^{m-1}\sqrt{C_{\mathrm{CY}}}
C_{\mathrm G}^{(m)}r^{1-m}.
\end{equation}
If \(m=2\), then
\[
1+|\ln d(x,z)|\leq 2+|\ln r|.
\]
By \eqref{equ:1 order} and Lemma~\ref{lem:bundle-green-from-heat}, we obtain
\begin{equation}\label{equ:first-derivative-case-two-dim2}
\max\bigl\{|\bar\partial_yG_{p,q}(x,y)|,
|\bar\partial_y^*G_{p,q}(x,y)|\bigr\}
\leq 4\sqrt{C_{\mathrm{CY}}}
C_{\mathrm G}^{(2)}\frac{1+|\ln r|}{r}.
\end{equation}
For \(m=2\), the right-hand side of
\eqref{equ:first-derivative-case-one-explicit} is bounded by
\[
6\sqrt{C_{\mathrm{CY}}}C_{\partial\mathrm G}
\frac{1+|\ln r|}{r}.
\]
Thus the first-order constant can be chosen as
\begin{equation}\label{equ:explicit-first-order-green-constant}
C_{\nabla\mathrm G}:=
\begin{cases}
2^{m-1}\sqrt{C_{\mathrm{CY}}}
\max\{3C_{\partial\mathrm G},C_{\mathrm G}^{(m)}\},&m>2,\\[1mm]
\sqrt{C_{\mathrm{CY}}}
\max\{6C_{\partial\mathrm G},
4C_{\mathrm G}^{(2)}\},&m=2.
\end{cases}
\end{equation}
The boundary-weighted first-order constant can be chosen as
\begin{equation}\label{equ:explicit-boundary-first-order-green-constant}
C_{\partial\nabla\mathrm G}
:=2^m\sqrt{C_{\mathrm{CY}}}
\max\{3C_{\partial^2\mathrm G},C_{\partial\mathrm G}\},
\end{equation}
by \eqref{equ:first-derivative-weighted-case-one} and
\eqref{equ:first-derivative-weighted-case-two}.
Consequently, the constants in Theorem~\ref{thm: green form complex boundary}
may be taken to be
\(C_{\mathrm G}^{(m)}\) in (i) when $m>2$ and
\(C_{\mathrm G}^{(2)}\) in (i) when $m=2$, \(C_{\partial\mathrm G}\) and
\(C_{\partial^2\mathrm G}\) in (ii), and
\(C_{\nabla\mathrm G}\) and \(C_{\partial\nabla\mathrm G}\) in (iii).
The single constant in the statement is the maximum of the applicable
constants in these three parts.
\end{proof}

\section{From Green Representations to Sobolev-Type Inequalities}\label{section:real inequality}

This section derives Corollaries~\ref{cor: complex laplace} and
\ref{cor: complex dbar} from the Green representation formulas and the
mapping properties of the associated volume and boundary kernels.

\subsection{Green representation formulas}\

We begin by fixing the kernel pairing and the boundary operator, and then
record the bundle-valued Green representation formulas. Put
\(
F^{p,s}:=\Lambda^{p,s}T^*M\otimes E
\)
for \(0\leq s\leq n\), and set $F^{p,-1}=F^{p,n+1}=0.$
Let
\(
K(x,y)\in
F_x^{p,q}\otimes(F_y^{p,s})^*
\simeq
\operatorname{Hom}(F_y^{p,s},F_x^{p,q})
\)
and \(v\in F_y^{p,s}\). Their \(F_x^{p,q}\)-valued partial pairing is
\begin{equation*}
\left\langle\!\left\langle v,K(x,y)
\right\rangle\!\right\rangle_y
:=
K(x,y)[v].
\end{equation*}
For the Green kernel, the derivatives in the dual variable are defined by
\begin{align*}
\bar\partial_yG_{p,q}:&=
\bigl((\bar\partial^*)^\vee\bigr)_yG_{p,q}
\in
\operatorname{Hom}
\bigl(F_y^{p,q+1},F_x^{p,q}\bigr),
\\
\bar\partial_y^*G_{p,q}
&:=
\bigl(\bar\partial^\vee\bigr)_yG_{p,q}
\in
\operatorname{Hom}
\bigl(F_y^{p,q-1},F_x^{p,q}\bigr).
\end{align*}
Here \(L^\vee\) denotes the formal transpose of a differential operator, characterized by
\[
\int_\Omega (L^\vee\eta)(u)\,dV
=
\int_\Omega \eta(Lu)\,dV
\]
for compactly supported smooth sections $\eta, u$.
Choose a boundary defining function
\(\rho\in C^\infty(\overline\Omega,\mathbb R)\) satisfying
\[
\Omega=\{\rho<0\},
\qquad
d\rho\neq0\quad\text{on }\partial\Omega.
\]
The pointwise adjoint of left exterior multiplication by
\(\bar\partial\rho\) is
\begin{equation*}
\bigl(\bar\partial\rho\wedge\bigr)^*
=
\iota_{(\bar\partial\rho)^{\sharp}},
\end{equation*}
where \((\bar\partial\rho)^{\sharp}\in T^{0,1}M\) is the Hermitian
dual of \(\bar\partial\rho\), and $\iota$ denotes interior multiplication.
Set
\[
\gamma_\rho^-u:=\iota_{(\bar\partial\rho)^{\sharp}}u,
\qquad
\gamma_\rho^+u:=\bar\partial\rho\wedge u,
\]
and define
{\small
\[
(\mathcal B_{p,q}u)(x):=\int_{\partial\Omega}\left(
\left\langle\!\left\langle\gamma_\rho^-u(y),\bar\partial_y^*G_{p,q}(x,y)\right\rangle\!\right\rangle_y-
\left\langle\!\left\langle\gamma_\rho^+u(y),\bar\partial_yG_{p,q}(x,y)\right\rangle\!\right\rangle_y
\right)\frac{dS(y)}{|d\rho(y)|}.
\]
}
\begin{thm}\label{thm:integral representations}
With the notation of Theorem
\ref{thm: green form complex boundary}, let \(\rho\) and $\mathcal B_{p,q}$ be as above.
Then the following identities hold.
\begin{itemize}
\item[(i)]
For every \(f\in C^2(\overline\Omega,F^{p,q})\) and \(x\in\Omega\),
\[
f(x)
=
\int_\Omega
\left\langle\!\left\langle
\square_{p,q}f(y),G_{p,q}(x,y)
\right\rangle\!\right\rangle_y\,dV(y)
+
(\mathcal B_{p,q}f)(x).
\]

\item[(ii)]
For every \(h\in C^1(\overline\Omega,F^{p,q})\) and \(x\in\Omega\),
\[
\begin{aligned}
h(x)
={}&\int_\Omega
\Bigl(
 \left\langle\!\left\langle
 \bar\partial h(y),\bar\partial_yG_{p,q}(x,y)
 \right\rangle\!\right\rangle_y+
 \left\langle\!\left\langle
 \bar\partial^*h(y),\bar\partial_y^*G_{p,q}(x,y)
 \right\rangle\!\right\rangle_y
\Bigr)\,dV(y)\\
&\quad +
(\mathcal B_{p,q}h)(x).
\end{aligned}
\]
\end{itemize}
\end{thm}

The identities in Theorem~\ref{thm:integral representations} follow directly
from repeated integration-by-parts. A closely related argument for differential forms appears in \cite[Lemma~6.2]{DHQ24}; we omit the analogous bundle-valued argument.

\subsection{Volume and boundary potential estimates}\ 

    We retain the constants \(r_0\) and \(C_0\) chosen in the proof of
    Theorem~\ref{thm:cheng yau general}; in particular,
    \eqref{equ:local-volume-single-constant} holds for every
    \(x\in\overline\Omega\) and \(0<\rho\leq r_0\). Set
    \begin{equation}\label{equ:global-volume-constant}
    C_{\mathrm{vol}}
    :=\max\left\{C_0,
    \frac{\operatorname{Vol}(\Omega)}{r_0^m}\right\}.
    \end{equation}
    Then, for every \(x\in\overline\Omega\) and
    \(0<\rho\leq D_\Omega\),
    \begin{equation}\label{equ:global-volume-bound-from-local}
    \operatorname{Vol}\bigl(B(x,\rho)\cap\Omega\bigr)
    \leq C_{\mathrm{vol}}\rho^m.
    \end{equation}

    For \(1<r\leq\infty\), let \(r'=r/(r-1)\), with the usual
    convention \(r'=1\) when \(r=\infty\), and set
    \[
    \Phi_2(\rho):=
    \begin{cases}
    \rho^{2-m},&m>2,\\
    1+|\ln\rho|,&m=2,
    \end{cases}
    \qquad
    \Phi_1(\rho):=
    \begin{cases}
    \rho^{1-m},&m>2,\\
    (1+|\ln\rho|)\rho^{-1},&m=2.
    \end{cases}
    \]
    Define the kernel norms
    \begin{equation}\label{equ:explicit-Riesz-Linfty-constants}
    \mathcal I_{\alpha,r}
    :=\sup_{x\in\Omega}
    \left(\int_\Omega\Phi_\alpha(d(x,y))^{r'}\,dV(y)\right)^{1/r'},
    \qquad \alpha\in\{1,2\},
    \end{equation}
    whenever \(r>m/\alpha\). The convention \(1/\infty=0\) will also
    be used below.

    \begin{lem}\label{lem:explicit-Riesz-kernel-bounds}
    Let \(\alpha\in\{1,2\}\), \(1<r\leq\infty\), and
    \(r>m/\alpha\). Put
    \[
    \sigma_{\alpha,r}
    :=m-(m-\alpha)r'>0.
    \]
    If \(m>2\), then
    \begin{equation}\label{equ:explicit-I-alpha-r-high-dim}
    \mathcal I_{\alpha,r}
    \leq
    \left(
    \frac{mC_{\mathrm{vol}}}{\sigma_{\alpha,r}}
    \right)^{1/r'}
    D_\Omega^{\alpha-m/r}.
    \end{equation}
    If \(m=2\), set
    \[
    A_\Omega:=1+|\ln D_\Omega|.
    \]
    Then
    \begin{equation}\label{equ:explicit-I-alpha-r-dim-two}
    \mathcal I_{\alpha,r}
    \leq
    (2C_{\mathrm{vol}})^{1/r'}
    D_\Omega^{\alpha-2/r}
    \sigma_{\alpha,r}^{-1/r'}
    \left(
    A_\Omega+
    \frac{\Gamma(r'+1)^{1/r'}}{\sigma_{\alpha,r}}
    \right).
    \end{equation}
    Thus all the geometric constants in these bounds are expressed in
    terms of the constants \(m,D_\Omega,r_0,C_0\) and
    \(\operatorname{Vol}(\Omega)\) already occurring in Section~\ref{sec:green form complex estimate}.
    \end{lem}
    \begin{proof}
    For \(x\in\Omega\), set
    \[
    V_x(t):=\operatorname{Vol}\bigl(B(x,t)\cap\Omega\bigr),
    \qquad 0<t\leq D_\Omega.
    \]
    If $\Psi$ is nonnegative and decreasing and
    $t^m\Psi(t)\to0$ as $t\downarrow0$, Stieltjes integration by parts and
    \eqref{equ:global-volume-bound-from-local} give
    \begin{equation}\label{equ:layer-cake-kernel-bound}
    \int_\Omega\Psi(d(x,y))\,dV(y)
    \leq mC_{\mathrm{vol}}
    \int_0^{D_\Omega}\Psi(t)t^{m-1}\,dt.
    \end{equation}

    Suppose first that \(m>2\). The layer-cake formula and
    \eqref{equ:layer-cake-kernel-bound} give
    \begin{align*}
    \int_\Omega d(x,y)^{-(m-\alpha)r'}\,dV(y)
    &\leq
    mC_{\mathrm{vol}}
    \int_0^{D_\Omega}
    t^{m-1-(m-\alpha)r'}\,dt=\frac{mC_{\mathrm{vol}}}{\sigma_{\alpha,r}}
    D_\Omega^{\sigma_{\alpha,r}},
    \end{align*}
    which proves
    \eqref{equ:explicit-I-alpha-r-high-dim}.

    Now suppose that \(m=2\). For \(0<t\leq D_\Omega\), write
    \(t=D_\Omega s\), where \(0<s\leq1\). The triangle inequality gives
    \[
    1+|\ln t|
    \leq 1+|\ln D_\Omega|+|\ln s|
    =A_\Omega-\ln s.
    \]
    Applying \eqref{equ:layer-cake-kernel-bound} to the decreasing majorant
    below gives
    \begin{align*}
    \int_\Omega\Phi_\alpha(d(x,y))^{r'}\,dV(y)
    &\leq 2C_{\mathrm{vol}}
    \int_0^{D_\Omega}
    t^{1-(2-\alpha)r'}
    \left(A_\Omega-\ln\frac{t}{D_\Omega}\right)^{r'}dt\\
    &=2C_{\mathrm{vol}}
    D_\Omega^{2-(2-\alpha)r'}
    \int_0^1s^{\sigma_{\alpha,r}-1}
    (A_\Omega-\ln s)^{r'}\,ds.
    \end{align*}
    With \(u=-\ln s\), the last integral becomes
    \[
    \int_0^\infty e^{-\sigma_{\alpha,r}u}
    (A_\Omega+u)^{r'}\,du.
    \]
    Minkowski's inequality for the measure
    $e^{-\sigma_{\alpha,r}u}\,du$ yields
    \begin{align*}
    &\left(
    \int_0^\infty e^{-\sigma_{\alpha,r}u}
    (A_\Omega+u)^{r'}\,du
    \right)^{1/r'}\leq
    \sigma_{\alpha,r}^{-1/r'}
    \left(
    A_\Omega+
    \frac{\Gamma(r'+1)^{1/r'}}{\sigma_{\alpha,r}}
    \right).
    \end{align*}
    Taking the \(r'\)-th root and the supremum over \(x\) proves
    \eqref{equ:explicit-I-alpha-r-dim-two}.
    \end{proof}
     
    The following lemma records the volume-potential estimate from
    \cite{DHQ24} used below. We include the argument because the nonendpoint bounds follow
    from H\"older's inequality and Fubini's theorem, whereas the critical bound
    uses the Hardy--Littlewood--Sobolev inequality.

    \begin{lem}\label{lem:DHQ-volume-mapping}
    For \(\alpha\in\{1,2\}\), define
    \[
    (T_\alpha h)(x)
    :=\int_\Omega \Phi_\alpha(d(x,y))|h(y)|\,dV(y).
    \]
    If \((a,b)\) is \((m,\alpha)\)-admissible, then there is a constant 
    $C:=C(\overline\Omega,g,\alpha,a,b)>0$ such that 
    \[
    \|T_\alpha h\|_{L^b(\Omega)}
    \leq C_{\alpha,a,b}\|h\|_{L^a(\Omega)}.
    \]
    \end{lem}
    
    \begin{proof}
Set
\[
K_\alpha(x,y):=\Phi_\alpha(d(x,y)),
\qquad x,y\in\Omega.
\]
For \(1\leq s<\infty\), define
\[
M_{\alpha,s}
:=
\max\!\left\{
\begin{aligned}
&\sup_{x\in\Omega}
 \left(\int_\Omega K_\alpha(x,y)^s\,dV(y)\right)^{1/s},\\
&\sup_{y\in\Omega}
 \left(\int_\Omega K_\alpha(x,y)^s\,dV(x)\right)^{1/s}
\end{aligned}
\right\}.
\]
Since \(d(x,y)=d(y,x)\), the two suprema are equal.

We first record the kernel-integrability ranges. If \(\alpha<m\), then
Lemma~\ref{lem:explicit-Riesz-kernel-bounds} implies
\begin{equation}\label{equ:kernel-uniform-Ls}
M_{\alpha,s}<\infty
\qquad
\text{for every }
1\leq s<\frac{m}{m-\alpha}.
\end{equation}
Indeed, when \(s>1\), set
\[
r:=\frac{s}{s-1}.
\]
Then \(r'=s\), and
\[
r>\frac{m}{\alpha}
\quad\Longleftrightarrow\quad
s<\frac{m}{m-\alpha}.
\]
Thus \eqref{equ:kernel-uniform-Ls} follows from
Lemma~\ref{lem:explicit-Riesz-kernel-bounds}. The case \(s=1\)
corresponds to \(r=\infty\).

When \((m,\alpha)=(2,2)\), one has
\[
K_2(x,y)=1+|\log d(x,y)|,
\]
and the same lemma gives
\begin{equation}\label{equ:log-kernel-all-finite-moments}
M_{2,s}<\infty
\qquad\text{for every }1\leq s<\infty.
\end{equation}

We next recall the generalized Young inequality in the precise form needed
here: Suppose \(1\leq s<\infty\) and \(M_{\alpha,s}<\infty\). Then
\begin{equation}\label{equ:generalized-Young-relation}
\|T_\alpha h\|_{L^b(\Omega)}
\leq M_{\alpha,s}\|h\|_{L^a(\Omega)}
\end{equation}
whenever
\begin{equation}\label{equ:Young-exponent-relation}
1\leq a\leq b\leq\infty,
\qquad
1+\frac1b=\frac1a+\frac1s.
\end{equation}

In particular, taking \(s=1\) shows that
\begin{equation}\label{equ:Talpha-La-La}
\|T_\alpha h\|_{L^a(\Omega)}
\leq M_{\alpha,1}\|h\|_{L^a(\Omega)}
\qquad (1\leq a\leq\infty).
\end{equation}
Since \(\Omega\) has finite volume,
\[
\|u\|_{L^b(\Omega)}
\leq
\operatorname{Vol}(\Omega)^{1/b-1/c}
\|u\|_{L^c(\Omega)}
\qquad (1\leq b\leq c\leq\infty).
\]
Thus, in proving the lemma, it suffices to consider \(b\geq a\).

We now verify every admissible range.

\medskip
\noindent
\emph{Case 1: \(\alpha<m\) and \(a=1\).}
In this case admissibility means
\[
1\leq b<\frac{m}{m-\alpha}.
\]
Take \(s=b\). Then \(M_{\alpha,b}<\infty\) by
\eqref{equ:kernel-uniform-Ls}, and
\[
1+\frac1b=1+\frac1s.
\]
The generalized Young inequality therefore gives
\[
\|T_\alpha h\|_{L^b(\Omega)}
\leq M_{\alpha,b}\|h\|_{L^1(\Omega)}.
\]

\medskip
\noindent
\emph{Case 2: \(1<a<m/\alpha\).}
Put
\[
b_*:=\frac{ma}{m-\alpha a}.
\]
Suppose first that \(a\leq b<b_*.\)
Define \(s\) by
\begin{equation}\label{equ:s-from-a-b}
\frac1s:=1+\frac1b-\frac1a,
\end{equation}
then
\[
s<\frac{m}{m-\alpha}.
\]
Thus \(M_{\alpha,s}<\infty\), and the generalized Young inequality gives
\[
\|T_\alpha h\|_{L^b(\Omega)}
\leq M_{\alpha,s}\|h\|_{L^a(\Omega)}.
\]

It remains to consider \(b=b_*\). This endpoint is included in the
definition of admissibility only when
\[
(m,\alpha)\neq(2,1).
\]
Since \(\alpha<m\), this implies \(m>2\), and hence
\[
\Phi_\alpha(d(x,y))=d(x,y)^{\alpha-m}.
\]
Covering $\overline\Omega$ by finitely many coordinate charts, applying
the Euclidean Hardy--Littlewood--Sobolev inequality to the near-diagonal
part of the kernel, and estimating the smooth off-diagonal remainder gives
a constant $C:=C(\overline\Omega,g,\alpha,a)>0$ such that
\begin{equation}\label{equ:HLS-critical-Talpha}
\|T_\alpha h\|_{L^{ma/(m-\alpha a)}(\Omega)}
\leq C
\|h\|_{L^a(\Omega)}.
\end{equation}

When \((m,\alpha)=(2,1)\), the kernel is instead
\[
\Phi_1(d(x,y))
=
\frac{1+|\ln d(x,y)|}{d(x,y)}.
\]
It belongs uniformly to \(L^s\) only for \(s<2\). Consequently, the above
Young argument gives precisely
\[
b<\frac{2a}{2-a},
\]
and does not give the critical equality, in agreement with the definition
of \((2,1)\)-admissibility.

\medskip
\noindent
\emph{Case 3: \(a=m/\alpha\) and \(\alpha<m\).}
The admissible target exponent \(b\) is finite. By
\eqref{equ:Talpha-La-La}, it is enough to consider \(b\geq a\). Define
\(s\) by \eqref{equ:s-from-a-b}. Since \(1/a=\alpha/m\),
\[
\frac1s
=
1+\frac1b-\frac{\alpha}{m}
=
\frac{m-\alpha}{m}+\frac1b.
\]
Because \(b<\infty\),
\[
s<\frac{m}{m-\alpha}.
\]
Also \(b\geq a\) implies \(s\geq1\). Therefore
\(M_{\alpha,s}<\infty\), and the generalized Young inequality yields
\[
\|T_\alpha h\|_{L^b(\Omega)}
\leq C_{\alpha,a,b}\|h\|_{L^a(\Omega)}.
\]

\medskip
\noindent
\emph{Case 4: \(a>m/\alpha\) and \(\alpha<m\).}
Let \(a'=a/(a-1)\), with \(a'=1\) when \(a=\infty\). The condition
\(a>m/\alpha\) is equivalent to
\[
a'<\frac{m}{m-\alpha}.
\]
Thus \(M_{\alpha,a'}<\infty\). H\"older's inequality gives, uniformly in
\(x\in\Omega\),
\begin{align*}
T_\alpha h(x)
&\leq
\left(\int_\Omega
K_\alpha(x,y)^{a'}\,dV(y)\right)^{1/a'}
\|h\|_{L^a(\Omega)}\\
&\leq M_{\alpha,a'}\|h\|_{L^a(\Omega)}.
\end{align*}
Consequently,
\[
\|T_\alpha h\|_{L^\infty(\Omega)}
\leq M_{\alpha,a'}\|h\|_{L^a(\Omega)}.
\]
All smaller target exponents follow from the finite-volume embedding.

\medskip
\noindent
\emph{Case 5: \((m,\alpha)=(2,2)\).}
Here
\[
K_2(x,y)=1+|\log d(x,y)|
\]
has uniformly bounded \(L^s\)-norms for every finite \(s\), by
\eqref{equ:log-kernel-all-finite-moments}.

If \(a=1=m/\alpha\), then every finite \(b\) is admissible. Taking \(s=b\)
in the generalized Young inequality gives
\[
\|T_2h\|_{L^b(\Omega)}
\leq M_{2,b}\|h\|_{L^1(\Omega)}
\qquad (1\leq b<\infty).
\]
If \(a>1\), then \(a'<\infty\), so H\"older's inequality gives
\[
\|T_2h\|_{L^\infty(\Omega)}
\leq M_{2,a'}\|h\|_{L^a(\Omega)}.
\]
The estimates for all finite \(b\) again follow from the finite-volume
embedding.

Combining the five cases proves the assertion for every
\((m,\alpha)\)-admissible pair \((a,b)\). All constants involved depend
only on \((\overline\Omega,g,\alpha,a,b)\).
\end{proof}

    The boundary terms are governed by the intrinsic Poisson-scale kernel
    \[
    P(x,y):=\delta(x)d(x,y)^{-m},
    \qquad (x,y)\in\Omega\times\partial\Omega.
    \]
    The following mapping estimate supplies the boundary term in both Sobolev
    inequalities.

    \begin{lem}\label{lem:poisson-kernel-mapping}
    Let
    \[
    (T_Pf)(x):=\int_{\partial\Omega}P(x,y)|f(y)|\,dS(y).
    \]
    Then
    \[
    T_P:L^1(\partial\Omega)\longrightarrow
    L^{\frac m{m-1},\infty}(\Omega),
    \qquad
    T_P:L^\infty(\partial\Omega)\longrightarrow L^\infty(\Omega)
    \]
    are bounded, where $L^{m/(m-1),\infty}$ is the corresponding weak
    Lebesgue space. Moreover, if $(a,b)$ is $m$-boundary-admissible, then
    there is a constant $C=C(\overline\Omega,g,a,b)>0$ such that
    \[
    \|T_Pf\|_{L^b(\Omega)}
    \leq C\|f\|_{L^a(\partial\Omega)}.
    \]
    \end{lem}
    \begin{proof}
    Since \(\delta(x)\) is the distance from \(x\)
    to the whole boundary,
    \[
    \delta(x)\leq d(x,y),
    \qquad
    P(x,y)=\frac{\delta(x)}{d(x,y)^m}
    \leq d(x,y)^{1-m}.
    \]
    Therefore, \cite[Lemma~2.13]{DHQ24} gives
    \[
    T_P:L^1(\partial\Omega)\longrightarrow
    L^{\frac m{m-1},\infty}(\Omega).
    \]
   
    We next prove 
    $$T_P:L^\infty(\partial\Omega)\longrightarrow L^\infty(\Omega)$$
    is bounded. Fix
    \(x\in\Omega\) and put \(\tau:=\delta(x)>0\). Since
    \(d(x,y)\geq\tau\) for every \(y\in\partial\Omega\), the sets
    \[
    A_j(x):=\left\{y\in\partial\Omega:
    2^j\tau\leq d(x,y)<2^{j+1}\tau\right\},
    \qquad j=0,1,2,\ldots,
    \]
    cover \(\partial\Omega\). On \(A_j(x)\),
    \[
    P(x,y)\leq\tau(2^j\tau)^{-m}.
    \]
    On the other hand, there is a constant $C_S:=C_S(\overline\Omega,g)>0$ such that 
    \[
    \operatorname{Area}(A_j(x))
    \leq C_S(2^{j+1}\tau)^{m-1}.
    \]
    Consequently,
    \begin{align}
    \int_{\partial\Omega}P(x,y)\,dS(y)
    &\leq
    \sum_{j=0}^\infty
    \tau(2^j\tau)^{-m}
    C_S(2^{j+1}\tau)^{m-1}
    =2^mC_S. \label{equ:Poisson-Linfty-kernel}
    \end{align}
    It follows that
    \begin{equation}\label{equ:Poisson-infty-infty}
    \|T_Pg\|_{L^\infty(\Omega)}
    \leq2^mC_S\|g\|_{L^\infty(\partial\Omega)}.
    \end{equation}

    The operator $T_P$ is positive and sublinear. Marcinkiewicz
    interpolation between the preceding weak endpoint and
    \eqref{equ:Poisson-infty-infty} yields, for $1<a<\infty$,
    \[
    T_P:L^a(\partial\Omega)\longrightarrow
    L^{ma/(m-1)}(\Omega).
    \]
    Since $\Omega$ has finite volume, this also gives every
    $1\leq b\leq ma/(m-1)$. For $a=1$, the weak endpoint embeds into
    $L^b(\Omega)$ for every $b<m/(m-1)$; for $a=\infty$,
    \eqref{equ:Poisson-infty-infty} and the finite-volume embedding give all
    $1\leq b\leq\infty$. These are exactly the
    $m$-boundary-admissible pairs.
    \end{proof}
\subsection{Sobolev-type inequalities}\ 

We now derive Corollaries~\ref{cor: complex laplace} and
\ref{cor: complex dbar}. Along \(\partial\Omega\), the pointwise adjointness
of exterior and interior multiplication gives
\[
|\gamma_\rho^-u|+|\gamma_\rho^+u|
\leq C|d\rho|\,|u|.
\]
The boundary-weighted estimate in
Theorem~\ref{thm: green form complex boundary}(iii) therefore implies
\begin{equation}\label{equ:boundary-operator-pointwise-majorization}
|\mathcal B_{p,q}u(x)|
\leq C\int_{\partial\Omega}
\frac{\delta(x)}{d(x,y)^m}|u(y)|\,dS(y)
=C T_P(|u|)(x).
\end{equation}
The zero-order representation in
Theorem~\ref{thm:integral representations}(i) and part~(i) of
Theorem~\ref{thm: green form complex boundary} consequently give
\[
|f(x)|
\leq C T_2(|\square_{p,q}f|)(x)+C T_P(|f|_{\partial\Omega})(x).
\]
Similarly, the first-order representation gives
\[
|f(x)|
\leq C T_1(|\bar\partial f|)(x)
+C T_1(|\bar\partial^*f|)(x)
+C T_P(|f|_{\partial\Omega})(x).
\]
Applying Lemmas~\ref{lem:DHQ-volume-mapping} and
\ref{lem:poisson-kernel-mapping} term by term proves the two corollaries in
their stated admissible ranges.

\section{Quantitative Solvability of the Top-Degree
\texorpdfstring{$\bar\partial$}{d-bar}-Equation}
    \label{sec:L^p estimate}

We conclude by applying the first-order Green kernel estimate in top
antiholomorphic degree.
\begin{proof}[Proof of Theorem~\ref{thm:complex L^p estimate}]
Put \(m=2n\) and, off the diagonal, set
\[
K_{p,n}(x,y):=\bar\partial_x^*G_{p,n}(x,y).
\]
This is a fiber map from \(F_y^{p,n}\) to \(F_x^{p,n-1}\). More precisely,
the Hermitian symmetry \(G_{p,n}(x,y)=G_{p,n}(y,x)^*\) and the
dual-variable convention give the type-correct identity
\begin{equation}\label{equ:top-degree-kernel-symmetry}
K_{p,n}(x,y)
=
\left.
\bigl(\bar\partial_z^*G_{p,n}(y,z)\bigr)^*
\right|_{z=x}.
\end{equation}
Consequently, Theorem~\ref{thm: green form complex boundary}(iii) gives a
constant
\(C_K=C_K(\overline\Omega,g,E)>0\) such that
\begin{equation}\label{equ:top-degree-kernel-bound}
|K_{p,n}(x,y)|\leq C_K
\begin{cases}
d(x,y)^{1-m},&m>2,\\
(1+|\ln d(x,y)|)d(x,y)^{-1},&m=2.
\end{cases}
\end{equation}
For smooth $f$, \(K_{p,n}\) is the Schwartz kernel of
\(\bar\partial^*\mathcal G_{p,n}\), and hence
\[
(\mathcal S_{p,n}f)(x)
=\int_\Omega K_{p,n}(x,y)f(y)\,dV(y).
\]
Thus
\[
|\mathcal S_{p,n}f(x)|
\leq C_K T_1(|f|)(x),
\]
so Lemma~\ref{lem:DHQ-volume-mapping} gives the asserted \(L^r\)-to-
\(L^k\) estimate for every admissible pair.

We next prove the right-inverse identity for smooth data. Let
$f\in C^\infty(\overline\Omega,F^{p,n})$ and set
$v:=\mathcal G_{p,n}f$. Dirichlet regularity gives
$v\in C^\infty(\overline\Omega,F^{p,n})$, and
$\bar\partial v=0$ for degree reasons. Hence
\[
\bar\partial(\mathcal S_{p,n}f)
=\bar\partial\bar\partial^*v
=(\bar\partial\bar\partial^*+\bar\partial^*\bar\partial)v
=\square_{p,n}v=f.
\]

It remains to pass to general data. Suppose first that \(1\leq r<\infty\),
and choose
\(f_j\in C^\infty(\overline\Omega,F^{p,n})\) with
\(f_j\to f\) in \(L^r\). The mapping estimate just proved shows that
\(\mathcal S_{p,n}f_j\) converges in \(L^k\); this defines the unique bounded
extension of \(\mathcal S_{p,n}\) to \(L^r\). Since the distributional
operator
\[
\bar\partial:L^k(\Omega,F^{p,n-1})
\longrightarrow\mathcal D'(\Omega,F^{p,n})
\]
is continuous, passing to the limit in
\(\bar\partial\mathcal S_{p,n}f_j=f_j\) yields
\[
\bar\partial\mathcal S_{p,n}f=f
\qquad\text{in }\mathcal D'(\Omega,F^{p,n}).
\]

Finally, let \(r=\infty\). Since \(\Omega\) has finite volume,
\(L^\infty(\Omega)\subset L^2(\Omega)\); define
\(\mathcal S_{p,n}f\) by restricting the \(L^2\)-extension constructed
above. This extension agrees almost everywhere with the kernel integral,
because the latter defines the same bounded \(L^2\) operator on smooth
forms. The integrability of the right-hand side of
\eqref{equ:top-degree-kernel-bound} gives directly
\[
\|\mathcal S_{p,n}f\|_{L^\infty(\Omega)}
\leq C\|f\|_{L^\infty(\Omega)},
\]
and the estimates for finite \(k\) follow by finite-volume embedding. The
distributional right-inverse identity already follows from the \(r=2\)
case. This proves all assertions, including both infinite endpoints.
\end{proof}

\end{document}